\documentclass{iopart}

\usepackage{iopams}
\expandafter\let\csname equation*\endcsname\relax
\expandafter\let\csname endequation*\endcsname\relax

\usepackage{amssymb,amsmath,amsfonts,amscd,verbatim}
\usepackage{url,upgreek}

\newcommand\textopenone{\leavevmode\hbox{\small 1\kern-3.3pt\normalsize 1}}
\newcommand\openone{\mathalpha \mbox{\textopenone}}
\newcommand{\myone}{\openone}

\usepackage{tikz}
\usetikzlibrary{matrix}
\usetikzlibrary{positioning}
\usepackage{pgfplots}
\usepackage{latexsym, amsmath, xcolor, graphicx, amssymb, color}

\pgfplotsset{compat=1.6}
\usepackage{wrapfig}
\usepackage{float}
\usepackage{mathtools}
\usepackage[top=1.2in,bottom=1.2in,left=1.2in,right=1.2in]{geometry}
\usepackage{algorithm,algorithmic}
\usepackage{cleveref}

\newcommand{\R}{\mathbb{R}}
\newcommand{\B}[1]{\mathbf{#1}}
\newcommand{\mc}[1]{\mathcal{#1}}

\newcommand{\bmat}[1]{\begin{bmatrix}#1 \end{bmatrix}}

\input{mymacros.sty}

\begin{document}
\title[Bayesian Level Set Approach]{Bayesian Level Set Approach for Inverse Problems with Piecewise Constant Reconstructions}
\author{William Reese${}^1$, Arvind K.\ Saibaba${}^1$, Jonghyun Lee${}^2$}

\address{$^1$ Department of Mathematics, 
North Carolina State University, Raleigh NC 27695, USA}
\address{$^2$ Department of Civil and Environmental Engineering and Water Resources Research Center Honolulu Hawaii 96822, USA}
\begin{abstract}
    There are several challenges associated with inverse problems in which we seek to reconstruct a piecewise constant field, and which we model using multiple level sets. Adopting a Bayesian viewpoint, we impose prior distributions on both the level set functions that determine the piecewise constant regions as well as the parameters that determine their magnitudes. We develop a Gauss-Newton approach with a backtracking line search to efficiently compute the maximum a priori (MAP) estimate as a solution to the inverse problem. We use the Gauss-Newton Laplace approximation to construct a Gaussian approximation of the posterior distribution and use preconditioned Krylov subspace methods to sample from the resulting approximation. To visualize the uncertainty associated with the parameter reconstructions we compute the approximate posterior variance using a matrix-free Monte Carlo diagonal estimator, which we develop in this paper. We will demonstrate the benefits of our approach and solvers on synthetic test problems (photoacoustic and hydraulic tomography, respectively a linear and nonlinear inverse problem) as well as an application to X-ray imaging with real data.
\end{abstract}
\vspace{2pc}

\submitto{\IP}
\section{Introduction}
Inverse problems allow us to estimate physical quantities that may not be directly observable~\cite{bardsley2018computational,vogel2002computational,hansen2010discrete,kaipio2005statistical}. Two examples that we consider in this paper are X-ray tomography and hydraulic tomography.  X-ray tomography is a non invasive technique of to image the inside of a body. X-ray sources are placed around the object and then radiation passes through the sources to corresponding detectors. The attenuation of the X-rays is used to reconstruct the internal details of the body \cite{kaipio2005statistical}. In hydraulic tomography a series of pumping tests is performed with multiple pumping and observation wells to determine complex subsurface structures such as spatially distributed hydraulic conductivity of an aquifer \cite{lee2013}. In both of these examples the unknown quantity (an image or hydraulic conductivity field) may be represented as a piecewise constant spatial field while the negligible intra-variability can be ignored. Beyond these two examples, inverse problems have applications in many areas of science and engineering such as biomedicine, geosciences, and material science.

 In deterministic inverse problems, the unknown parameter field is uniquely determined such that the estimated field can reproduce the observations within a reasonable error level. In addition measurements may be sparse or limited and the recovery of the unknown parameter field from the measurements may be an ill-posed problem. Such is the case, for example, in limited angle X-ray tomography or subsurface structure characterization in hydraulic tomography. Uncertainty Quantification (UQ) allows us to address these different challenges. In deterministic inverse problems the ill-posedness is often alleviated by solving a regularized optimization problem in which the regularization term (e.g., Tikhonov regularization) often penalizes undesired behavior in the reconstructions. Alternatively, the Bayesian approach assumes the unknown field is random and characterizes the prior beliefs as a prior distribution. The prior distribution is combined with the likelihood distribution using the Bayes rule to produce the posterior distribution which serves as the solution to the Bayesian inverse problem. In this paper, we are primarily interested in inverse problems in which the unknown fields are assumed to be piecewise constant. For example in hydraulic tomography, one would like to identify a subsurface permeability field with discrete geological zones or facies, such as a regional aquifer formed with deposited layers of sediment, that are relatively homogeneous with each zone but with abrupt changes at the boundaries between zones. Widely used methods Tikhonov regularization or Gaussian-priors require a large amount of data and have high computational costs associated with the data size to provide sharply changing parameter estimation. This motivates advanced, scalable techniques to solve inverse problems with piecewise constant fields.

The literature on inverse problems with piecewise constant unknowns is vast and so we only mention a few key references. Perhaps the most common is the use of total variation (TV)~\cite{hansen2010discrete,vogel2002computational,kaipio2005statistical,bardsley2018computational}.  TV regularization is favored over Tikhonov regularization when the unknown parameter is discontinuous because it preserves sharp interfaces as opposed to smoothing them. A difficulty of using the TV regularization is the fact that it is non differentiable; however, there are many ways this can be handled (e.g., using smoothed TV approximations). The TV approach can also be used in the Bayesian setting by incorporating the TV regularization as a prior information (this can correspond to the Laplace prior). However, in \cite{siltanen2004}, the authors demonstrate that the conditional mean estimates for the TV prior are not edge preserving through fine discretizations of the model space. There are other approaches for promoting edge preserving behavior such as Besov priors~\cite{bui2015scalable}, level set approaches~\cite{blip,iglesias2016,dunlop2016map}, and Cauchy difference priors~\cite{MarkkanenRoininenHuttunenLasanen+2019+225+240}. In this paper, we develop a multiple level set based approach for enforcing piecewise constant reconstructions. The multiple level set approach utilizes a combination of level set functions, that determine the shapes, and constants, representing the magnitudes, to model the unknown parameter field. The advantage of level set functions is their ability to model the discontinuous parameter field in terms of a differentiable function. This changes the inversion problem from estimating a piecewise constant field into a few (possibly smooth) level set functions that determine the regions and the magnitudes of the field in those regions.

\paragraph{Overview of the main contributions.}  This paper proposes a Bayesian level set approach for solving inverse problems in which the unknown parameters are represented as piecewise constant fields. There are three main contributions of this paper:
\begin{enumerate}
    \item  We use  multiple level sets to enforce piecewise constant reconstructions. This approach allows the flexibility to incorporate prior information for the level set functions $\vPhi$ which control the shape and the parameters governing the magnitudes of the representations $\vc$. We discuss several possible choices for prior distributions. 
    \item Assuming a Gaussian distribution for the measurement errors and the prior distribution, we derive the posterior distribution, which, in general is non-Gaussian because of the nonlinear parameterization. We develop an inexact Gauss-Newton approach to obtain the maximum a posteriori (MAP) estimate. The algorithm is efficient in that it can take advantage of the adjoint-based method for computing derivatives and only incurs a small additional cost due to the level set parameterization. 
    \item To quantify the uncertainty in the reconstructed parameters, we approximate the posterior distribution using Laplace's approximation based on the Gauss-Newton Hessian; this is known as linearized Bayesian inference (see e.g.,~\cite{villa2021hippylib}). Using this approximate distribution we develop matrix-free (here and henceforth, by matrix-free, we mean that the approaches only rely on matrix-vector products and do not require us to form the operators explicitly) iterative methods to generate samples using the preconditioned Lanczos approach and estimate the posterior variance using two different methods. As with the MAP estimate, the posterior samples also have the piecewise constant structure. We also develop a matrix-free estimator for the diagonals of the inverse of a matrix (LanczosMC) that combines the Lanczos and Monte Carlo approaches, and {which may be of independent interest beyond this paper}. 
\end{enumerate} 
Finally, we demonstrate the benefits of the Bayesian level set approach and the computational benefits of the proposed algorithms on a range of applications such as Photoacoustic Tomography, X-ray Tomography using a real-world data set, and hydraulic tomography (which is a nonlinear inverse problem).

 \paragraph{Related approaches.} We position our work in relation to other works on level methods in inverse problems.  Surveys of level set methods in deterministic inverse problems are provided in~\cite{burger2005survey,dorn2006level}. Another deterministic approach that uses level sets is parametric level set approach (PaLS)~\cite{aghasi2011parametric}. In PaLS, the level set function(s) are represented a weighted sum of radial basis functions and the unknown parameters are the weights. The advantage of this approach over our multiple level set approach is the use of a few parameters (e.g., on the order of $10-100$) rather than a pixel or voxel-wise approach that can require millions of unknowns.  Our work is closely related to \cite{van2010multiple}, who also use a multiple level set approach and a Gauss-Newton type approach for solving the optimization problem.  However, quantifying the reconstruction uncertainty is still an active area of research. 
 
 Several authors have utilized level sets to perform Bayesian inversion for piecewise constant parameters.  Cardiff and Kitanidis~\cite{blip} demonstrate the use of level set functions for facies detection in a hydraulic tomography problem and provide uncertainty estimates for the facies' boundaries. This approach relies on computing of the velocity field for the level set functions to update the interfaces. In \cite{iglesias2016} the authors formulate a Bayesian level set inverse problem in terms of a single level set function, but assume that the number of levels and their magnitudes are known {\em a priori}. The authors formulate the Bayesian inverse problem in an infinite dimensional setting and utilize Monte Carlo Markov Chain (MCMC) methods to sample from the posterior distribution. This approach does not require computing the velocity field for the level set function evolution as it is updated through the sampling. In \cite{dunlop2017hierarchical} the authors extend the previous approach by incorporating length and amplitude scales as parameters in addition to the level set function in the prior.  In contrast to these approaches, which use MCMC techniques to explore the posterior distribution, our approach uses matrix-free techniques that are scalable to large problem sizes (both in the number of unknown parameters and measurements).

\paragraph{Overview of the paper.} We conclude this section with a brief overview of this paper. In Section~\ref{sec:back}, we review the basics of inverse problems and the level set approach. In Section~\ref{sec:bayeslevelset}, we present the Bayesian level set approach, discuss several options for prior information, and use an inexact Gauss-Newton approach for computing the maximum a posteriori (MAP) estimate. In Section~\ref{sec:uq}, we derive an expression for the approximate posterior distribution and develop matrix-free techniques to generate samples from this distribution and the estimate the posterior variance. In Section~\ref{sec:numex}, we apply this framework and the solvers we have developed to synthetic test problems from photoacoustic and hydraulic tomography, and real-world problem from X-ray tomography. In \ref{sec:diaginv}, we derive a matrix-free technique to estimate the diagonals of the inverse of a matrix and validate on several different matrices; we dissociate it from the rest of the paper since it may be of interest beyond this paper. Finally, we offer some concluding remarks and future directions in Section~\ref{sec:conc}.

\section{Background}\label{sec:back}
In this section, we give background information on the Bayesian inverse problem and level set technique for representing piecewise constant functions.

\subsection{Bayesian Inverse Problem}\label{ssec:bayes}
Let $\vf \in \R^{n} \rightarrow \R^{n_{\rm obs}}$ represent the forward operator or the parameter-to-observable map, let $\vm \in \R^n$ denote the unknown parameters to be recovered. If we assume that the measurement are corrupted by additive Gaussian noise $\boldsymbol\varepsilon$, then the following equation relates the measurements $\B{d} \in \R^{n_{\rm obs}},$ to the unknown parameters 
\begin{equation}\label{eqn:meas}
\vd = \vf(\vm) + \boldsymbol\varepsilon \;\;\;\; \boldsymbol\varepsilon \sim \mathcal{N}(\B{0},\boldsymbol\Gamma_\text{noise}).
\end{equation}
Here, we represent the measurement error $\boldsymbol\varepsilon$ as a Gaussian distribution with zero mean, and a positive definite covariance matrix $\mgamma_\text{noise}$. The inverse problem can be stated as follows: given measurements $\vd$ use \eqref{eqn:meas} to recover the unknown parameters $\vm$.

Fundamental limitations in data acquisition, means that the recovery of spatial parameters $\vm$ is mathematically an {\em ill-posed} problem: there may be no solution, the solution may not be unique, or it may depend sensitively on the data $\vd$. The Bayesian approach specifies a prior probability distribution $\pi_\mathrm{pr}(\vm)$ that describes the prior information, or ``belief,'' about the parameter $\vm$. Then Bayes' rule is applied to derive the posterior probability distribution, defined as the distribution of  the parameter $\vm$ conditioned on the data $\B{d}$ and is represented by $\pi_\text{post}(\vm\mid\vd)$. The posterior distribution takes the form
\begin{equation}\label{eqn:post}
\pi_\text{post}(\vm \mid \vd) \propto  \>  \pi_\text{like}(\vd\mid\vm) \pi_\text{pr}(\vm) \propto \exp\left( -\frac12\|\vd - \vf(\vm)\|_{\boldsymbol\Gamma_\text{noise}^{-1}}^2\right) \pi_\text{pr}(\vm).
\end{equation}

A standard approach is to compute the so-called {\em maximum a posteriori} (MAP) estimate. That is, the MAP  point $\vm_\text{MAP}$ maximizes the posterior distribution $\pi_\text{post}(\vm \mid \vd)$. If the prior distribution is of the form $\pi_\text{prior}(\vm) \propto\exp(-\mc{R}(\vm))$, the MAP estimate is the solution to the variational problem
\begin{equation}\label{eqn:det}
\min_{\vm \in \R^n} \> \frac12\|\vd - \vf(\vm)\|_{\boldsymbol{\Gamma}_\text{noise}^{-1}}^2 + \mathcal{R}(\vm). 
\end{equation}
The MAP estimate is easy to interpret as the ``most likely'' value of the parameter $\vm$, given the data $\vd$. The MAP estimate can be computed using optimization methods and the specific choice of the solver (either a general purpose or a special purpose) depends on the choice of the prior distribution~\cite{stuart2010inverse}.

\subsection{Level Set Representation}
To explain the idea behind the level set method, consider a piecewise constant function $m$ that only takes two positive values $c_p$ and $c_b$, in the domain $\Omega \subset \real^d$; here the dimension takes values $d=1,2,3$. The function can be represented mathematically as follows. Let $A \subset \Omega$ be a closed subset and let $\partial \Omega$ denote its boundary. Let $\chi_A$ denote an indicator function, which takes two values $1$ and $0$, depending on whether the point $\vr\in A$ or outside of $A$, respectively. Then, the function can be written in terms of the indicator function as 
\[ m(\vr) = c_p\chi_A(\vr) + c_b[1-\chi_A(\vr)] \qquad \forall \vr \in \Omega .\]

\begin{figure}[!ht]
    \centering
    \includegraphics[scale=0.3]{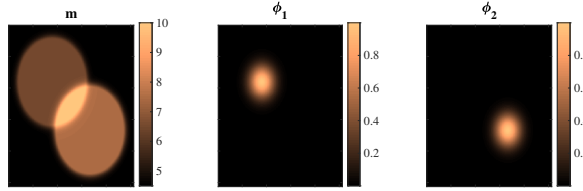}
    \caption{Left: piecewise constant field with 4 regions, Right: level set functions $\phi_1$ and $\phi_2$ that are used to compute $m$.}
    \label{fig:2circlelvlset}
\end{figure}

Then, in order to completely specify the function $m$, we only need knowledge of the region $A$. The level set approach is shape-based and represents the boundary $\partial A$ as the level set of a function $\phi:\Omega \rightarrow \real$. Specifically, if we take the zero level set, then $\phi(\vr)$ is related to $A$ and $\partial A$ as 
\[  \begin{array}{ll}   \phi(\vr) > 0 &  \forall\vr \in A \\ \phi(\vr) = 0 & \forall\vr \in \partial A \\ \phi(\vr) < 0 & \forall\vr \in \Omega \backslash A.\end{array} \]
Alternatively, we can represent the function $m(\vr)$ in terms of the Heaviside function $H(x) = \frac{1}{2}(1+\text{sign}(x))$  as 
\[  m(\vr) = c_p H(\phi(\vr)) + c_b[1- H(\phi(\vr))].\]
By this formulation, identifying the unknown region $A$, is tantamount to estimating the continuous function $\phi(\vr)$ and the coefficients $c_p$ and $c_b$. To represent $m(\vr)$ when it has 4 regions we can use 2 level set functions and write 
\begin{align*}
m(\vr) &= c_1 (1-H(\vphi_1))(1-H(\vphi_2)) + c_2 H(\vphi_1)(1-H(\vphi_2)) \\
				&+ c_3 (1-H(\vphi_1))H(\vphi_2) + c_4H(\vphi_1)H(\vphi_2)
\end{align*}
We illustrate this representation for two level set functions in Figure~\ref{fig:2circlelvlset}.

We now explain how to extend this approach to recover piecewise constant functions with more than $2$ regions. This follows the discussion in~\cite{van2010multiple}. Suppose that there are $\nr$ regions to be recovered. We choose the number of level set functions $\nls$ to be such that 
\begin{equation} 2^{\nls-1} < \nr \leq 2^{\nls}.\label{eqn:nls} \end{equation}
Denote the corresponding level set functions $\phi_i(\vr)$ for $1 \leq i \leq\nls$.  Then, the unknown function $f(\vr)$ can be expressed in terms of the unknown level set functions as 
\begin{equation}\label{eq:levelsetparam}
m(\vr) = \sum_{i_1,\dots,i_{\nls} = \{0,1\}} c_{i_1,\dots,i_{\nls}} \widehat{H}_{i_1}(\phi_1(\vr)) \cdots \widehat{H}_{i_{\nls}}(\phi_{\nls}(\vr)),
\end{equation}
where $c_{i_1,\dots,i_{\nls}}$ are the magnitudes of the piecewise constant regions, and summation involves $2^{\nls}$ terms. Furthermore, the modified level set functions $\widehat{H}_i(\cdot)$ have the form 
\[ \widehat{H}_i(x) = \left\{ \begin{array}{ll} H(x) & i = 0 \\ 1 - H(x) & i = 1.  \end{array}\right.\]
In this formulation, recovering the unknown function $m(\vr)$ now involves recovering the level set functions $\phi_i(\vr)$ and the coefficients $c_{i_1,\dots,i_{\nls}}$. In practice, the Heaviside function $H(x) = \frac{1}{2}(1+\text{sign}(x))$ is not differentiable; to address this issue we replace it with the mollified Heaviside function $\hsmooth$ as suggested in~\cite{aghasi2011parametric}:
\begin{equation}\label{eq:mollifiedHeav}
    \hsmooth(x) =  \left\{\begin{array}{ll}
    1 ,  & x > \epsilon ,  \\
    \frac12 + \frac{x}{2\epsilon} + \frac{1}{2\epsilon}\sin(\frac{x\pi}{\epsilon}), & |x| \leq \epsilon ,\\
        0 ,  & x < -\epsilon ,
    \end{array}\right.
\end{equation}
where $\epsilon$ is a positive parameter used to control the resolution of $\hsmooth$ at the point $x=0$. This choice of the Heaviside function ensures sparsity in the Jacobian which results in lower computational costs.
\begin{figure}[!ht]
    \centering
    \includegraphics[scale=0.5]{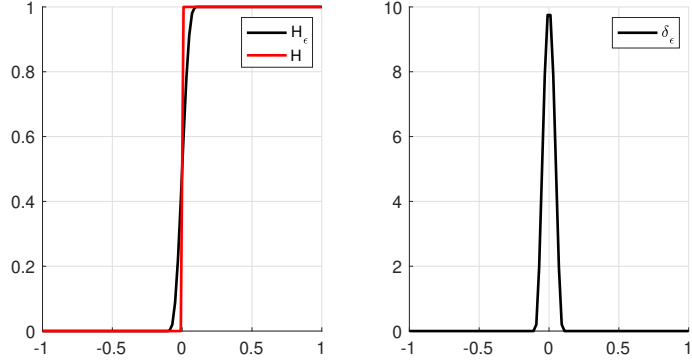}
    \caption{(left)  Heaviside function $H(x)$ and its approximation $\hsmooth(x)$, (right) derivative of Heaviside function $\delta_\epsilon(x)$.}
    \label{fig:heaviside}
\end{figure}
\paragraph{Discretization} The unknown coefficients $c_{i_1,\dots,i_{\nls}}$ can be rearranged in column major order into a vector $\vc$ of the size $2^{\nls}$. 
The domain $\Omega$ is discretized into $N$ points, and the vectors $\vphi^{(1)},\dots,\vphi^{(\nls)}$ represent the discretized versions with respect to some basis.  Concatenating the vectors $\vphi^{(i)}$ into the vector 
$$
\vPhi = \begin{bmatrix}
\vphi^{(1)} \\
\vdots \\
\vphi^{(\nls)}
\end{bmatrix} \in \real^{n{\nls}}
$$
and utilizing \eqref{eq:levelsetparam} we can write the unknown vector $\vm$ as 
\begin{equation}\label{eq:disc_parameter_form}
    \vm = \vm(\vPhi,\vc)
\end{equation}
i.e., $\vm$ is a vector valued function $\vm: \real^{(n\nls)}\times\real^{2^{\nls}}\longrightarrow\rn$. Now that we have a way to present piecewise constant images we now address the problem of estimating the parameters $\vPhi$ and $\vc$, which describe the image, from a set of measurements $\vd$. To this end, we adopt the Bayesian approach described in~\cite{kaipio2005statistical,stuart2010inverse}.

\section{Bayesian Level Set Formulation and MAP Estimate}\label{sec:bayeslevelset}

With the parameterization of $\vm$ as defined in \eqref{eq:disc_parameter_form}, we can rewrite the measurement equation as 
\begin{equation}\label{eq:meas}
    \vd = \vf(\vm(\vPhi,\vc)) + \vnoise.
\end{equation}
As in Section~\ref{ssec:bayes}, we assume that $\vnoise$ is Gaussian with a known mean and covariance $\vnoise \sim \pN(\vzero,\mnoise)$. This equation defines the likelihood $$\pi(\vd|\vPhi,\vc) \propto \exp\left(-\frac12\noisenorm{ \vf(\vm(\vPhi,\vc))-\vd} \right).$$
 We assume that our unknown parameters $\vPhi$ and $\vc$ are independent random vectors and our goal is to estimate these unknown parameters, given measurements $\vd$.  The posterior distribution combines the likelihood and the prior information through Bayes' theorem. However, the specific form of the posterior distribution depends on the specific choice of the prior distribution. We briefly discuss various possible choices for the prior distributions. 
\subsection{Choice of prior distributions}
Since our unknowns consist of the coefficients $\vc$ and the level set functions $\vPhi$ it is reasonable to assume that they are independently distributed random vectors. Thus we will discuss possible choices of prior for $\vc$ and $\vPhi$ separately. 

\subsubsection{Prior for the level set functions {$\vPhi$}}
Following the level set formulation it is apparent that we want each of the level set functions $\vphi_i$ to be ``smooth'' for $1 \leq i \leq \nls$. One approach to enforce smoothness in prior distributions is to use Gaussian distribution, which we adopt in this paper; it takes the form 
\begin{equation}\label{eq:phiprior}
    \vphi \sim \pN(\vmeanphi,\regphi^{-2}\mgammaphi)
\end{equation}
where $\regphi^{-2}$ is the precision and $\mgammaphi \in \rnn$ is the covariance matrix and $\vmeanphi$ is the mean. 
To assign a prior for $\vPhi$ we assume that each level set function is independently distributed and each level set function has the distribution given by \eqref{eq:phiprior}. In other words we assume $\vphi_i\sim \pN(\vmeanphi,\regPhi^{-2}\mgammaphi)$,  are independent for $1\leq i \leq \nls.$ Thus, we can write the prior for $\vPhi$ as 
\begin{equation}\label{eq:Phiprior}
\vPhi \sim \pN(\vmeanPhi,\regPhi^{-2}\mgammaPhi)
\end{equation}
where the mean $\vmeanPhi = \myone \otimes \vmeanphi$  (here \myone= $[{1,\dots,1}]^T \in \R^{\nls}$), 
and the covariance matrix $\mgammaPhi$ takes the form $\mgammaPhi = \mi_{\nls} \otimes \mgammaphi.$
The symbol $\otimes$ represents the Kronecker product. 
We assume that each level set function is independent and Gaussian with mean $\vmeanphi \in \rn$ and covariance $\mgammaphi \in \rnn$, so that we assume the same regularization parameter $\regPhi = \regphi$ for each level set. This reduces the number of hyperparameters to be estimated.

There are several possible choices for the covariance matrix $\mgammaphi$. If we assume that the underlying field $\phi$ is a Gaussian random field, then the covariance matrix $\mgammaphi$ can be chosen based on a kernel function. There are several choices here, e.g., Mat\'ern class, $\gamma$-exponential, etc~\cite{rasmussen2003gaussian}. Another approach is to choose priors from the Whittle-Mat\'ern class as was done in~\cite{lassi2014whittle}. Our choice of covariance matrix $\mgammaphi$ will depend on the test problem.
Yet another choice, which we adopt in this paper, is to use the Gauss-Markov Random field (GMRF) to specify the covariance matrix~\cite{bardsley2018computational,GMRFbook}. We choose the precision matrices (inverse of the covariance matrix $\mgammaphi$) as the discretization of a partial differential operator. This has close connections with the Gaussian random field approach~\cite{lindgren2011explicit}.

\subsubsection{Prior for level set magnitudes {$\vc$}}
The elements of $\vc$ represent the magnitudes of the values taken by $\vm$. When $\vm$ represents a piecewise constant image, the elements $c_i$, $1\leq i \leq 2^{\nls}$ assign intensity of pixels belonging to different regions. Therefore, imposing a prior distribution for $\vc$ controls the distribution of the intensity of pixels. We briefly mention a few choices for prior distributions.  A convenient representation for the prior information for the intensities is the multivariate Gaussian distribution, i.e.,
\begin{equation}\label{eq:cprior}
    \vc \sim \pN(\vmeanc,\mgammac)
\end{equation}
where $\vmeanc \in \real^{2^{\nls}}$ and $\mgammac \in \real^{{2^{\nls}}\times {2^{\nls}}}$ are the mean and covariance respectively. This is the approach we adopt in this paper.

Alternatively, if one requires that $\vc$ have positive elements, a lognormal or gamma prior may be more appropriate. In the lognormal case, we assume that $\vy=\log\vc \sim \pN(\vmean_c, \mgammac)$ (where the log is computed elementwise). Finally, if we assume the pixel intensities are independent and identically distributed according to the gamma distribution with its shape and (inverse) scale parameters $\alpha > 0$ and $\beta > 0$. That is, 
$c_i \sim \text{Gamma}(\alpha,\beta), 1 \leq i \leq \nls.$ This implies that the pdf of the prior distribution satisfies  
\begin{equation*}\label{eq:c_gamma}
    \pi(\vc) \propto  \prod_{i=1}^{2^{\nls}} \Bigg (c_i^{\alpha-1}\text{exp}(- \beta c_i)\Bigg).
\end{equation*}

\subsection{Efficient solver for computing the  MAP estimate}\label{ssec:gn}
For the rest of the section, we assume that the prior distribution for $\vPhi$ and $\vc$ are Gaussian. That is, we assume the following model
\begin{equation}
    \vPhi \sim  \mathcal{N}(\vmeanphi, \lambda^{-2}_{\vPhi}\mgammaphi) \qquad    \vc \sim \mathcal{N}(\vmeanc, \lambda^{-2}_{\vc}\mgammac).
\end{equation}

With this choice of prior distribution, the posterior distribution takes the form 
\begin{equation}\label{eq:post}
\post(\vPhi,\vc\mid\vd) \propto \bigexp{-\frac12 \noisenorm{ \vf(\vm(\vPhi,\vc))-\vd}-\frac{\regPhi^2}{2} \Phinorm{\vPhi-\vmeanPhi} -\frac{\lambda^{2}_{\vc}}{2} \cnorm{\vc-\vmeanc}}.
\end{equation}
 It is worth noting that due to the parameterization, the posterior distribution is non-Gaussian, even though the measurement noise and the prior distributions are Gaussian. This is due to the nonlinearity in the forward operator as well as the nonlinearity in the level set representation $\vm(\vPhi,\vc)$.

The MAP point is a maximizer of the posterior distribution; alternatively, it can be obtained as the minimizer of the negative logarithm of the posterior distribution. More specifically, the MAP point can be found by solving the optimization problem
\begin{equation}\label{eq:opt}
\min_{\vPhi,\vc} \> \left[-\log\post(\vPhi,\vc\mid\vd) \right] = \min_{\vPhi,\vc}\textrm{F}(\vPhi,\vc)
\end{equation}
where 
\begin{equation}\label{eq:functional}
    \textrm{F}(\vPhi,\vc)= \frac12 \noisenorm{ \vf(\vm(\vPhi,\vc))-\vd}+\frac{\lambda^2_{\vPhi}}{2}\Phinorm{\vPhi-\vmeanPhi}+\frac{\lambda^2_{\vc}}{2}\cnorm{\vc-\vmeanc}.
\end{equation}
Since~\eqref{eq:opt} is a nonlinear least squares problem we solve this efficiently using an inexact Gauss-Newton approach.

For simplicity, we denote the concatenated vector $\vx = \begin{bmatrix} \vPhi^T & \vc^T \end{bmatrix}^T.$ Using the chain rule, the Jacobian of $\vf$ with respect to $\vx$ is
\begin{equation}\label{eq:jacobian}
\mj = \p{\vf}{\vm}\p{\vm}{\vx}.
\end{equation}
The term $\p{\vf}{\vm}$ is the part of the Jacobian that comes from the forward operator and is problem dependent. It is important to note that this matrix need not be computed explicitly but we only need to be able to compute its action, i.e., compute matrix-vector products with this matrix. The Jacobian matrix of $\vm$ with respect to $\vx$ has the following form 
\begin{equation}\label{eq:dmdphic}
\p{\vm}{\vx} = \begin{bmatrix}\p{\vm}{\vphi_1} &  \dots &\p{\vm}{\vphi_{\nls}} &\p{\vm}{\vc}\end{bmatrix} 
\end{equation}
where the matrices
\begin{equation}\label{eq:dmdphi}
\p{\vm}{\vphi_i} \in \mathbb{R}^{n\times n} \qquad  1 \leq i \leq \nls \qquad \text{and}  \quad \p{\vm}{\vc} \in \mathbb{R}^{n\times 2^{\nls}}.
\end{equation}
 The details are given in~\cite{van2010multiple}. 
 
\begin{figure}[H]
    \centering
    \includegraphics[trim=0cm 1.8cm 0cm 1cm, clip=true, width=0.9\textwidth]{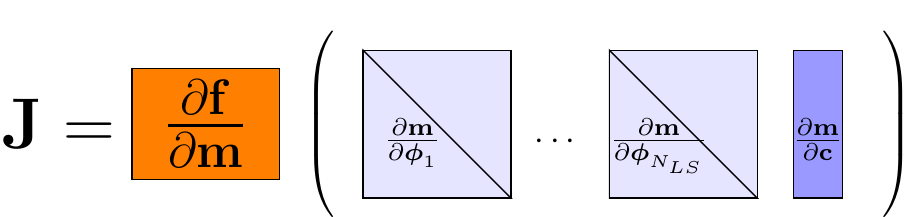}
   \caption{Visualization of the structure of the Jacobian.}
    \label{fig:jacobian}
\end{figure}

It is worth noting that the matrix $\p{\vm}{\vx}$ has a lot of sparsity that can be exploited for computational advantages. The matrices $\p{\vm}{\vphi_i}$ are diagonal and can be stored in sparse format. There may be additional sparsity in the diagonal entries and in ${\p{\vm}{\vc}}$ depending on the support of the derivative of the mollified Heaviside function. In Figure~\ref{fig:jacobian}, we illustrate the structure of the Jacobian. 

To compute the MAP estimate, we solve the nonlinear least squares problem \eqref{eq:opt} with an inexact Gauss-Newton approach~\cite[Chapter 10.2]{nocedal1999}. Given the current estimate $\vx_k$, the gradient of the objective function is given by
\begin{equation}\label{eq:gradient}
\vg_k = \mj^T_k\mnoise^{-1}\left(\vf(\vm_k)-\vd \right)+ \mprior^{-1}(\vx_k - \vprior)
\end{equation}
where $\vm_k = \vm(\vPhi_k,\vc_k)$ and
\begin{equation}\label{eq:priorMatrix}
\vprior = \begin{bmatrix} \vmeanPhi \\ \vmeanc\end{bmatrix}	\qquad \mprior =
	\begin{pmatrix}
	\lambda^{-2}_{\vPhi}\mgammaPhi & \\ & \lambda^{-2}_{\vc}\mgammac
	\end{pmatrix}.
\end{equation}

The Gauss-Newton step updates the current iterate as $\vx_{k+1} = \vx_k + \alpha_k\delta \vx_k$ where $\alpha_k$ is a step size and the search direction $\delta\vx_k$ is obtained by solving the linear system of equations
\begin{equation}\label{eq:gaussNewtonstep}
    \mh_k
	\delta\vx_k
	 = -\vg_k \qquad \text{where} \qquad \mh_k = \mj_k^T\mnoise^{-1}\mj_k + \mprior^{-1}.
\end{equation}
The matrix $\mh_k$ is known as the Gauss-Newton Hessian at step $k$.  The system of equations is solved using the conjugate gradient (CG) iterative method. The stopping criterion for the CG iterations is $\| \mh_k\delta \vx_k + \vg_k\|_2 \leq \eta_k \| \vg_k\|_2$ where the residual tolerance is taken to be a decreasing sequence of the form
\begin{equation}\label{eq:forcingSequence}
    \eta_k = \min\left(0.5,\sqrt{\frac{\|\vg_k\|_2}{\|\vg_0\|_2}}\right).
\end{equation}

The step length $\alpha_k$ is chosen using a backtracking line search procedure in order to satisfy the {strong Wolfe Conditions} with  parameters $c_1= 10^{-4}$, $c_2=0.9$~\cite[Section 3.1]{nocedal1999}. To stop the inexact Gauss-Newton iterations we use a combination of two different stopping criteria. The first criterion is based on the Morozov discrepancy principle~\cite[Chapter 7]{vogel2002computational} and terminates the iterations if the data misfit falls below a specified amount, whereas the second criterion is based on the first order optimality condition  and terminates the iterations if the gradient is smaller than a user-defined tolerance.

\subsection{Computational cost}\label{ssec:compcosts}
 To assess computational costs of each Gauss-Newton iteration we consider storage of matrices as well as the number of matrix-vector products required. In terms of storage, the additional storage required is $\nls$ diagonal $n \times n$ matrices ${\p{\vm}{\vphi_i}}$ with sparse diagonals, and one dense $n \times 2^{\nls}$ matrix ${\p{\vm}{\vc}}$; the total cost of storage is $n(\nls + 2^{\nls})$. Note that the matrix $\mj_k$ need not be formed explicitly.

Let $N_\text{CG}$ represent the number of CG iterations in one Gauss-Newton iteration. Let $T_\text{der}$ be the cost of applying the derivative ${\p{\vf}{\vm}}$ or its adjoint, and let $T_\text{prior}$ be the cost of applying the precision matrix $\mprior^{-1}$. Then, the cost of applying the Jacobian $\mj_k$ or its adjoint is 
\[T_\text{jac} = T_\text{der} + \mc{O}(n(\nls + 2^{\nls}))  \> \text{flops}.\] Since each Gauss-Newton iteration requires solving the system \eqref{eq:gaussNewtonstep} we need to compute the gradient $\vg_k$ and matrix vector products with $\mh_k$. The gradient involves a matrix-vector product with the adjoint of $\mj_k$ and $\mprior^{-1}$. Each matrix vector product with $\mh_k$ requires a matrix vector product with $\mprior^{-1}$, one with $\mj_k$ and its adjoint. Therefore, the total cost of a Gauss-Newton step (not including step length computation) is 
\[ T_\text{GN} = N_\text{CG}(2T_\text{jac} + T_\text{prior}) + T_\text{jac} + T_\text{prior} \> \text{flops}.\]
An important point to make here is that the use of the level set approach has only a marginal increase in the computational cost per iteration and the memory footprint. Furthermore, the proposed approach can be easily incorporated into standard solvers for inverse problems, since it reutilizes the ability to apply Jacobians and their adjoints.   

\section{Linearized Bayesian Inference}\label{sec:uq}

A common approach for uncertainty quantification is to generate independent samples from the posterior distribution and use these samples along with Monte Carlo estimates to generate statistics of the quantities of interest. Due to the nonlinearity of the Heaviside function, the posterior distribution is, in general, non-Gaussian even if the measurement operator $\vf(\vm) = \ma\vm$ is linear and the prior distributions are Gaussian. This implies that the problem of generating samples from the posterior distribution can be challenging. The prevalent approach for sampling is some version of Markov Chain Monte Carlo approach. However, the performance of MCMC approaches suffers in high dimensional spaces and good proposal distributions are necessary for successful implementation of MCMC methods. In this paper, we settle for linearized Bayesian inference (following the terminology in~\cite{buithanh2013computational,villa2021hippylib}) and approximate the posterior distribution using the Laplace's approximation (with the Gauss-Newton Hessian instead of the full Hessian) . We develop matrix-free techniques for generating samples from the approximate distribution and to estimate the posterior variance. 

A word on the notation: We let $\vx = \bmat{\vPhi^T & \vc^T}^T$ as the concatenated unknown vector and consider the mapping $\vm(\vx)$ as the mapping from the unknown vector to the pixel space. With this notation, the posterior distribution can be written as $\pi_\mathrm{post}(\vx|\vd)$.

\subsection{Sampling from the posterior distribution}\label{ssec:sample} We first derive an expression for the approximate posterior distribution by using the Laplace's approximation to the posterior distribution. This gives us a Gaussian distribution centered about the MAP estimate $\vmean_{\vx} = [\vPhi_\text{MAP}^T, \vc_\text{MAP}^T]^T$ and with the covariance $\mgamma_{\vx}$, given by the expression
\[ \mgamma_{\vx}= \> (\mjhat^T\mnoise^{-1}\mjhat + \mprior^{-1})^{-1},\]
where $\mjhat$ is the Jacobian of $\vf(\cdot)$ computed at $\vmean_{\vx}$. In other words, the covariance matrix is taken to be the inverse of the Gauss-Newton Hessian  computed at the MAP point, i.e., $\mgamma_{\vx} = \mh_\text{GN}^{-1}$. Note that we are using the Gauss-Newton Hessian rather than the full Newton Hessian since the former is guaranteed to be positive definite in the context of our problem. We denote this approximate posterior distribution $\vx|\vd \sim \mathcal{N}(\vmean_\vx,\mgamma_{\vx})$ with the resulting pdf denoted $\widehat\pi_{\vx}(\vx|\vd)$.

An important point to note here is that because of the size of the problems involved, forming the Gauss-Newton Hessian explicitly is not possible. As a result, explicit computation and storage of $\mgamma_{\vx}$ (or an appropriate factorization for it) is not feasible for large-scale problems since it is likely to be dense; for the same reason, computing the factor $\ms$ is also computationally infeasible. Therefore, to generate samples we use the preconditioned Lanczos approach outlined in~\cite{chow2014preconditioned}. The key advantage of this approach is that it only requires matrix-vector products with $\mgamma_{\vx}$ or its inverse.

To explain this approach, consider a factorization of the form $\mgamma_\vx = \ms\ms^T$. Note that we will not compute $\ms$ explicitly but use it only for illustration purpose. Given this factorization, we can generate samples from this distribution using the formula
\begin{equation}\label{eqn:xsample}
    \vx = \vmean_\vx + \ms\boldsymbol\eta \qquad \boldsymbol\eta \sim \mathcal{N}(\mathbf{0},\mi)
\end{equation}
Suppose we have the preconditioner $\mg$ such that $\mg\mg^T \approx \mgamma_{\vx}^{-1}$ and let
\[ \ms = \mg^{-T}(\mg^{-1}\mh_\text{GN}\mg^{-T})^{-1/2},\]
where $\mh_\text{GN} = \mgamma_{\vx}^{-1}$ is the Gauss-Newton Hessian evaluated at the MAP point. It is easy to verify that 
\[ \ms\ms^T = \mg^{-T}(\mg^{-1}\mh_\text{GN}\mg^{-T})^{-1/2}(\mg^{-1}\mh_\text{GN}\mg^{-T})^{-1/2}\mg^{-1} = \mgamma_{\vx}. \]
With this definition of $\ms$, the use of \eqref{eqn:xsample} requires the action of the matrix $(\mg^{-1}\mh_\text{GN}\mg^{-T})^{-1/2}$ on the vector $\boldsymbol\eta$.  This is accomplished using the Lanczos algorithm and it only requires matrix-vector products with the matrix $\mg^{-1}\mh_\text{GN}\mg^{-T}$.    See~\cite[Section 2]{chow2014preconditioned}; the details are omitted here. The goal of the preconditioner here is to reduce the number of iterations required by the Lanczos algorithm. In practice, throughout this paper, we choose the preconditioner $\mg$ to be the Cholesky factorization of $\mprior^{-1} = \mg\mg^T$. Such a factorization is easy to compute since $\mprior^{-1}$ is block-diagonal with sparse diagonal blocks.

Let $N_\text{Lanc}$ be the number of Lanczos iterations and let $T_\text{prec}$ be the cost of applying a preconditioner. Similar to the analysis in Section~\ref{ssec:compcosts}, the cost of generating a sample is 
\[ T_\text{sample} = N_\text{Lanc}(2T_\text{jac} + T_\text{prior} + 2T_\text{prec}) + \mathcal{O}(n\nls N_\text{Lanc})\> \text{flops}. \]
Since we use full reorthogonalization in the Lanczos process, there is an additional cost of $ \mathcal{O}((n\nls) N_\text{Lanc}^2)$ flops.

\subsection{Approximating the posterior variance}\label{ssec:postvar} 
A useful way to visualize the uncertainty associated with the reconstruction is to plot the posterior variance. Since the posterior distribution is non-Gaussian, computing the conditional variance is not straightforward. Here, we derive two methods for approximating the posterior variance based on the Gauss-Newton-Laplace approximation. First, we derive an approximate posterior distribution for $\vm$, conditioned on the data $\vd$.   Since the mapping $\vm(\vx)$ is differentiable, using a first-order Taylor expansion about the posterior mean $\vmean_{\vx}$
\[ \vm(\vx) \approx \vm(\vmean_{\vx}) + \mjm(\vx - \vmean_{\vx}), \]
where $\mjm = \p{\vm}{\vx}$ is the Jacobian of $\vm$ with respect to $\vx$ computed at the MAP estimate. Let $\vmean_{\vm} = \vm(\vmean_{\vx})-  \mjm\vmean_{\vx}$, then the approximate posterior distribution is 
\[ \vm|\vd \sim \mathcal{N}(\vmean_{\vm},\mgamma_{\vm})\qquad  \mgamma_{\vm} = \mjm \mgamma_{\vx}\mjm^\top, \]
with the corresponding pdf denoted $\widehat{\pi}_{\vm}(\vm|\vd)$. 
An approximation of the posterior variance can, therefore, be obtained from the diagonals of the matrix $\mgamma_{\vm}$. However, forming $\mgamma_{\vm}$ is computationally infeasible for the size of the problems, since we cannot compute $\mgamma_{\vx}$ explicitly. We now propose two different matrix-free methods for estimating the posterior variance. An advantage of both approaches is that they are embarrassingly parallel across the random vectors.

\subsubsection{Method 1: Using LanczosMC diagonal estimators} 
In the appendix, we propose a new method LanczosMC for estimating the diagonals of the inverse of a matrix. This method is a hybrid method which combines the advantages of the Lanczos method for estimating the diagonals and the Monte Carlo diagonal estimator.

A detailed algorithm, explanation of computational costs, as well as examples on test matrices will be give in \ref{sec:diaginv}. Here we only give a high-level overview. Our approach  first approximates $\mh_\text{GN}^{-1}$ using a low-rank approximation computed using the preconditioned Lanczos approach so that $\mh_\text{GN}^{-1} \approx \mw_k\mw_k^T$. Recall that $\mgamma_{\vx}=\mh_\text{GN}^{-1}$, which we can write as 
$$
\mh_\text{GN}^{-1} = \mw_k\mw_k^T +  \left(\mh_\text{GN}^{-1} - \mw_k\mw_k^T\right).
$$
The definition of $\mgamma_{\vm}$ implies that 
\begin{equation}\label{eq:diag_eq}
    \diag(\mgamma_{\vm}) = \diag(\mjm\mw_k\mw_k^T\mjm^T) + \diag(\mgamma_{\vm} - \mjm\mw_k\mw_k^T\mjm^T).
\end{equation}
The first summand is easy to calculate because it is a low-rank outer product. The second summand can be estimated with a Monte Carlo diagonal estimator, similar to the one described in~\cite{bekas2007estimator}. To estimate $\diag(\mgamma_{\vm} - \mjm\mw_k\mw_k^T\mjm^T)$, we draw $N$ independent random vectors $\{\vz_j\}_{j=1}^N$ with the Rademacher distribution (i.e., entries $\pm 1$ with equal probability). We then estimate 
\[\diag(\mgamma_{\vm} - \mjm\mw_k\mw_k^T\mjm^T) \approx \frac{1}{N} \sum_{j=1}^N\vy_j \odot \vz_j,  \]
where $\odot$ denotes the elementwise product and 
\[ \vy_j = \mjm (\mh_\text{GN}^{-1} (\mjm^T \vz_j)) - \mjm\mw_k\mw_k^T(\mjm^T \vz_j),  \qquad 1 \leq j \leq N.\]
Here to compute $\mh_\text{GN}^{-1} (\mjm^T \vz_j)$ we use the CG approach similar to that used in Section~\ref{ssec:gn}. Note that the only source of error in \eqref{eq:diag_eq} is in the second summand because the diagonal of a low-rank outer product in the first summand can be computed analytically. In \ref{sec:diaginv}, we show through a series of numerical experiments that the LanczosMC approach has lower variance compared to directly applying the Monte Carlo diagonal estimator to $\mgamma_{\vm}$. The cost of solving a linear system is comparable to that of computing the Gauss-Newton step; we denote it by $T_{\rm solve}$. The matrix $\mj_{\vm}$ is very sparse, so the cost of applying $\mj_{\vm}$ is negligible compared to the other costs. Therefore, the cost of using this approach is 
\[ T_\text{var1} = k(2T_\text{jac} + 2T_\text{prec} + T_\text{prior}) + NT_{\rm solve}+  \mc{O}(n\nls(Nk + k^2))\> \text{flops}.\]

\subsubsection{Method 2: Using approximate posterior samples} An alternative method to estimate the posterior variance is to use the samples generated from the approximate posterior distribution with pdf $\widehat\pi_{\vx}(\vx|\vd)$. Suppose we have independent samples $\{\vx_j\}_{j=1}^N$ from the distribution with the pdf $\widehat\pi_{\vx}(\vx|\vd)$ which can be computed using the methods proposed in Section~\ref{ssec:sample}, then define the sample average
$\widehat{\vm} = \frac1N \sum_{j=1}^N\vm(\vx_j)$. Using the multivariate delta method~\cite[Section 1.8, Theorem 8.22]{lehmann2006theory}, it can be shown that  
\[ \sqrt{N}(\widehat{\vm}-\vm(\vmean_{\vx})) \longrightarrow \mc{N}(\mathbf{0}, \mjm\mgamma_{\vx}\mjm^T),\]
as $N\rightarrow\infty$. Here the convergence is in the sense of distributions. As before, $\mj_{\vm}$ is the Jacobian computed at the MAP estimate. Therefore, we can approximate the posterior variance as the diagonals of the sample covariance operator 
\[ \diag(\mgamma_{\vm}) \approx \frac{1}{N-1}\diag\left(  \sum_{j=1}^N(\vm(\vx_j) - \widehat{\vm} )(\vm(\vx_j) - \widehat{\vm} ))^\top \right).    \]
The cost of this approach is simply $T_\text{var2} = NT_\text{sample}$ flops, where $T_\text{sample}$ is the cost of computing a sample using the approach described in Section~\ref{ssec:sample}.

\section{Numerical Experiments}\label{sec:numex}
In this section, we demonstrate the Bayesian level set approach and the performance of our solvers on a set of test problems. 
\subsection{Description of the test problems and parameters}
We describe the various test problems and images that we use in this paper. The first test application (Section~\ref{ssec:pat}) comes from Photoacoustic Tomography. The linear forward operator matrix $\ma$ is of size $23168 \times 16384$ and is generated using the \texttt{PRspherical} function from the IR Tools toolbox~\cite{gazzola2019ir}. We use two different test images from IR tools: the first one is called `Three Phases', which has three regions,  and the other is called `Grains'. 
The second application comes from X-ray tomography and uses an open-source tomographic dataset corresponding to a walnut image~\cite{hamalainen2015tomographic}. A brief summary of the problem sizes and the number of level set functions used are summarized below:
\begin{center}
\begin{tabular}{|c||c|c|c|}
\hline
 \textbf{Image}  & 
   \textbf{Three Phases} & \textbf{Grains} & \textbf{Walnut} \\
  
   \hline
   Application & \multicolumn{2}{|c|}{Photoacoustic} & {X-ray}  \\
      \hline
   $n_{\rm obs}$ &   $23168$ &  $23168$ & $9840$ \\
    \hline
   $n $ & $128 \times 128$ & $128\times 128$ & $82 \times 82$ \\
    \hline
    $\nls$  & $2$ & $3$ & $2$ \\\hline
\end{tabular}
\end{center}

\paragraph{Level set parameters} In the experiments, we choose the scale parameter $\epsilon=10^{-2}$. We experimented with different values of $\epsilon$ and observed that with a smaller $\epsilon$, the images were a lot sharper but consequently the optimization problem became harder. The number of level sets used depended on the test image used; in applications, a domain expert may have some prior knowledge on the maximum number of regions $N_R$ possible, and  one can use $\nls \sim \log_2 (N_R)$ using the formula in~\eqref{eqn:nls}. In the X-ray dataset, the true image represents a walnut and even though we do not know the exact number of regions, our approach does a reasonable job estimating the image.

\paragraph{Prior information} In all the experiments, the prior distributions for $\vPhi$ is a Gaussian distribution with zero mean and  covariance defined by $\mgammaPhi^{-1} = \mi_{\nls} \otimes \regPhi^2\left(\alpha\ml + \gamma\mi\right)^2$, where $\ml$ is the finite difference discretization of $(-\Delta)$ with homogeneous Neumann boundary conditions. The choice of the precision matrix for each individual level set is related to~\cite{bui2015scalable,bardsley2018computational}. The parameters $\alpha$ and $\gamma$ control the correlation length and smoothness of the samples drawn from the prior and we take $\alpha= 0.01$ and $\gamma = 0.1$. 
We also take the prior information for $\vc$ as a Gaussian distribution with zero mean and covariance matrix $\mgammac=\lambda_{\vc}^2\mi$. We discuss how to estimate $\regPhi^2$ in the experiments below. 

\subsection{Photoacoustic Tomography}\label{ssec:pat}
The first set of numerical experiments is based on a test problem from the  \verb|IRTools| package~\cite{gazzola2019ir} that models photoacoustic tomography (PAT). In all the experiments in this subsection, we used a regular grid of size $n = 128^2$. The matrix $\ma$ has the size $23168\times 16384$. We use two different test images: `Three Phases' and `Grains'. The noise level is set to be $2\%$ for both the test images. For the `Three Phases' problem, we take $\nls =2$, which is capable of handling $2^{\nls} = 4$ regions. On the other hand for the `Grains' test problem, which has multiple regions, we pick $\nls = 3$ which is capable of $8$ distinct regions. For both problems take $\lambda_{\vc}^2 = 10 \lambda_{\vPhi}^2$.  
The parameter $\lambda_{\vPhi}^2$ for each test problem is determined using an L-curve type approach. In this approach, we solve the inverse problem for several values of $\lambda_{\vPhi}^2$ and for each solution, plot the data misfit $\noisenorm{\vd-\ma\vm(\vx)}$ and the regularization term (without the scaling factor $\lambda_{\vPhi}^2$). We also plot the estimate of noise,  and we choose $\lambda_{\vPhi}^2$ at the first instance at which the solution is smaller than the estimate of the noise.  An alternative approach would be to use the continuation approach described in~\cite{haber2000optimization}.

\begin{figure}[!ht]
    \centering
    \includegraphics[scale=0.5]{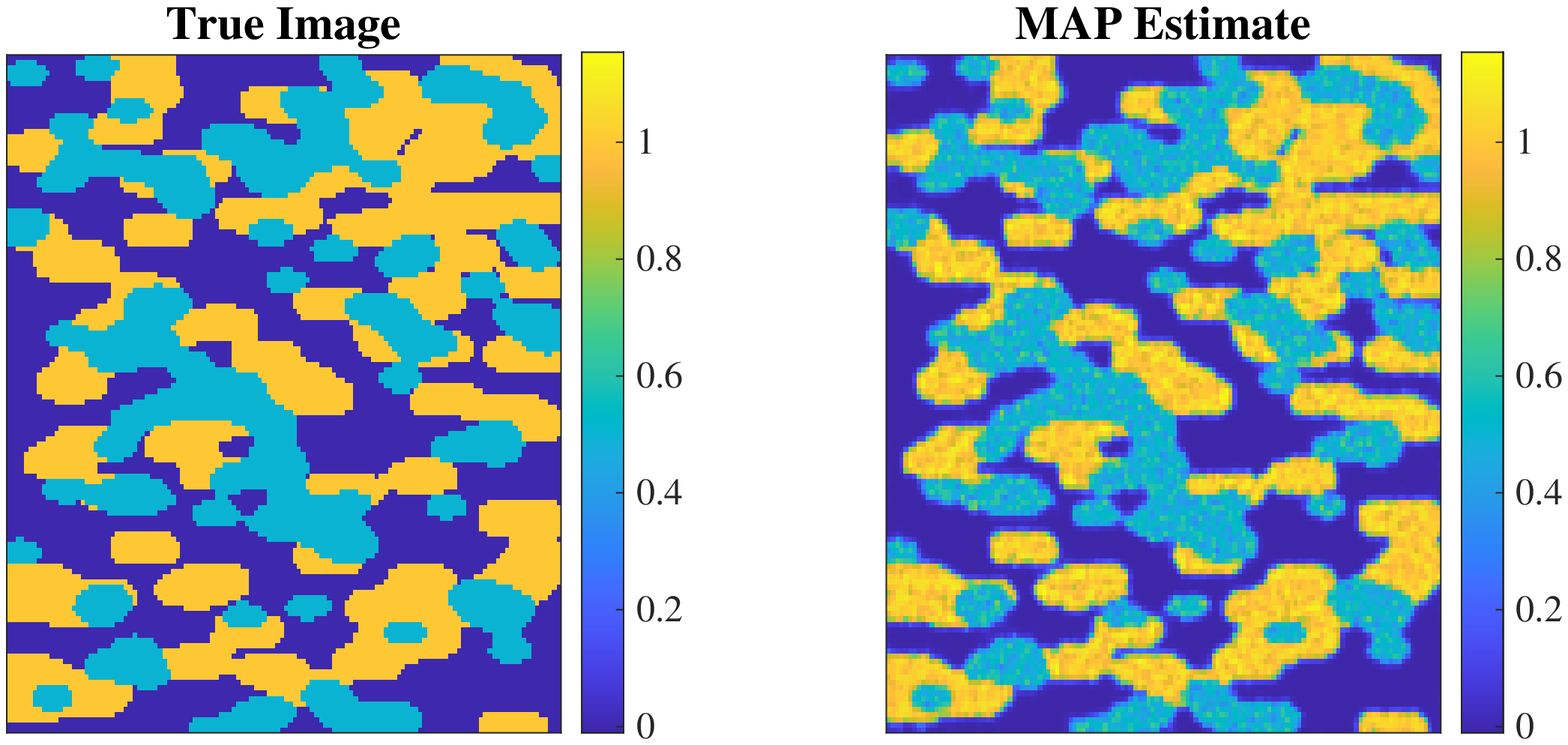} \\
    \includegraphics[scale=0.5]{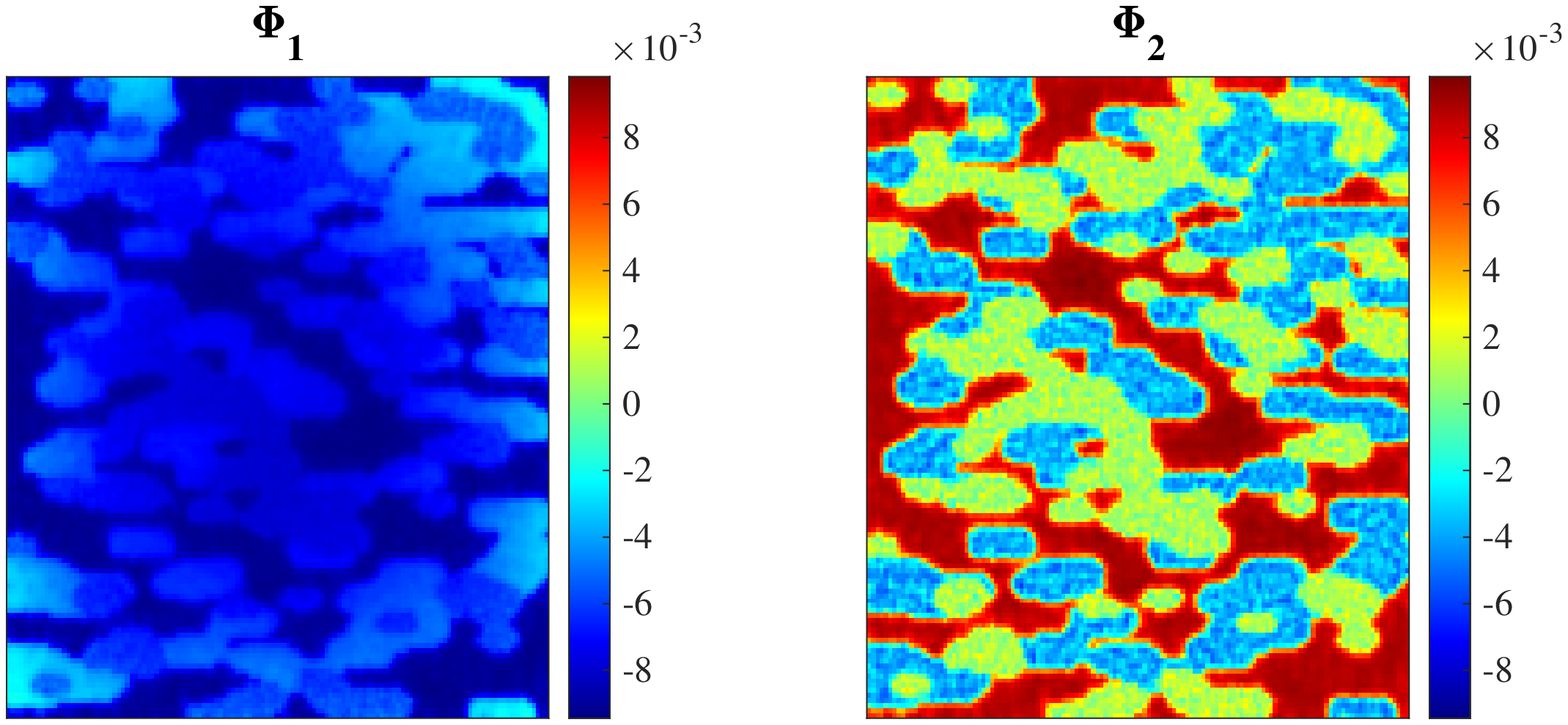}
    \caption{These images correspond to the `Three Phases' test problem. In the top row, we plot the true image (left) and the MAP estimate (right). In the bottom row, we plot the level set functions $\vphi_1$ and $\vphi_2$ that determine the MAP estimate. }
    \label{fig:threephasesmap}
\end{figure}

\begin{figure}[!ht]
    \centering
    \includegraphics[scale=0.25]{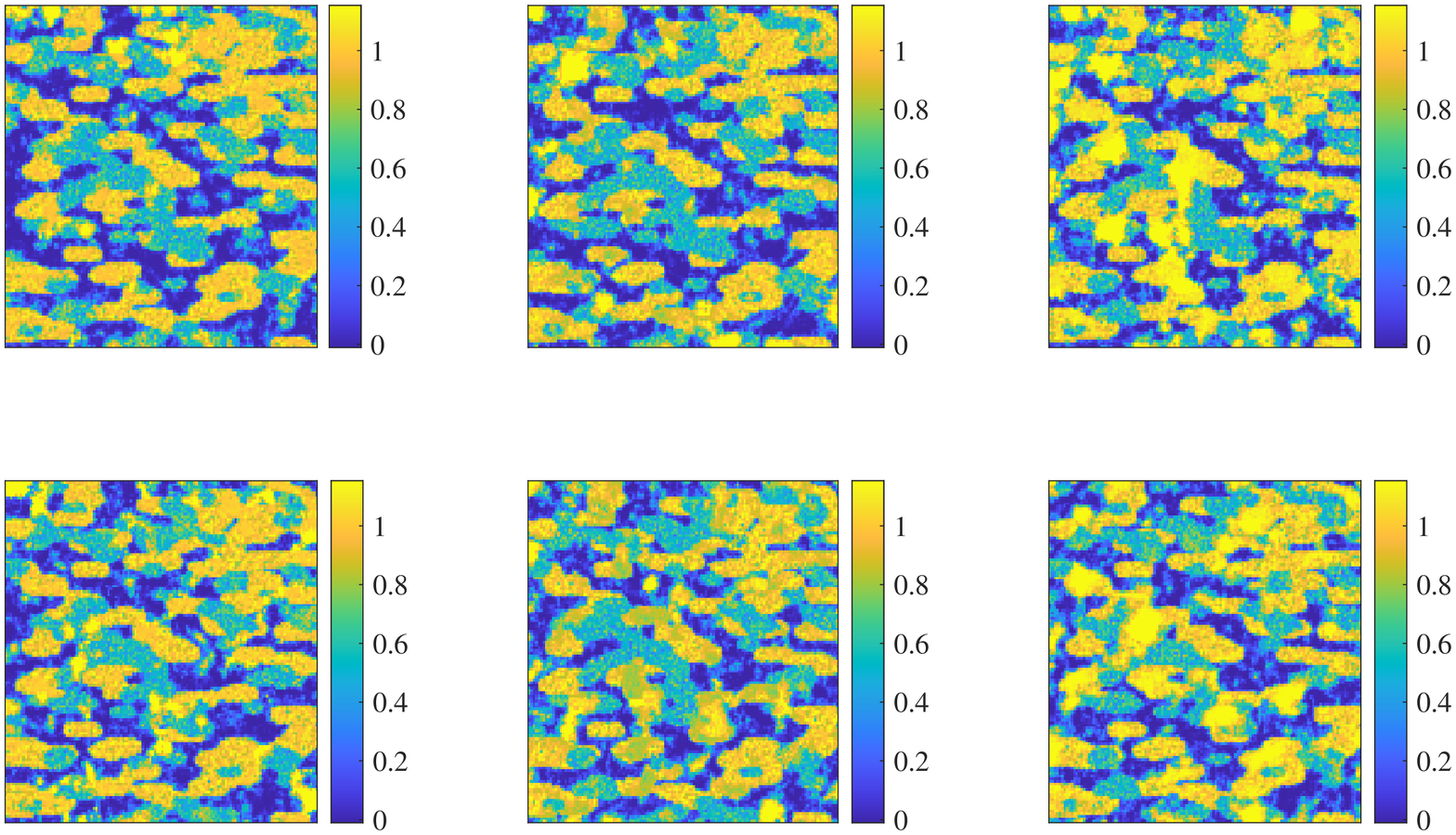}
    \includegraphics[scale=0.26]{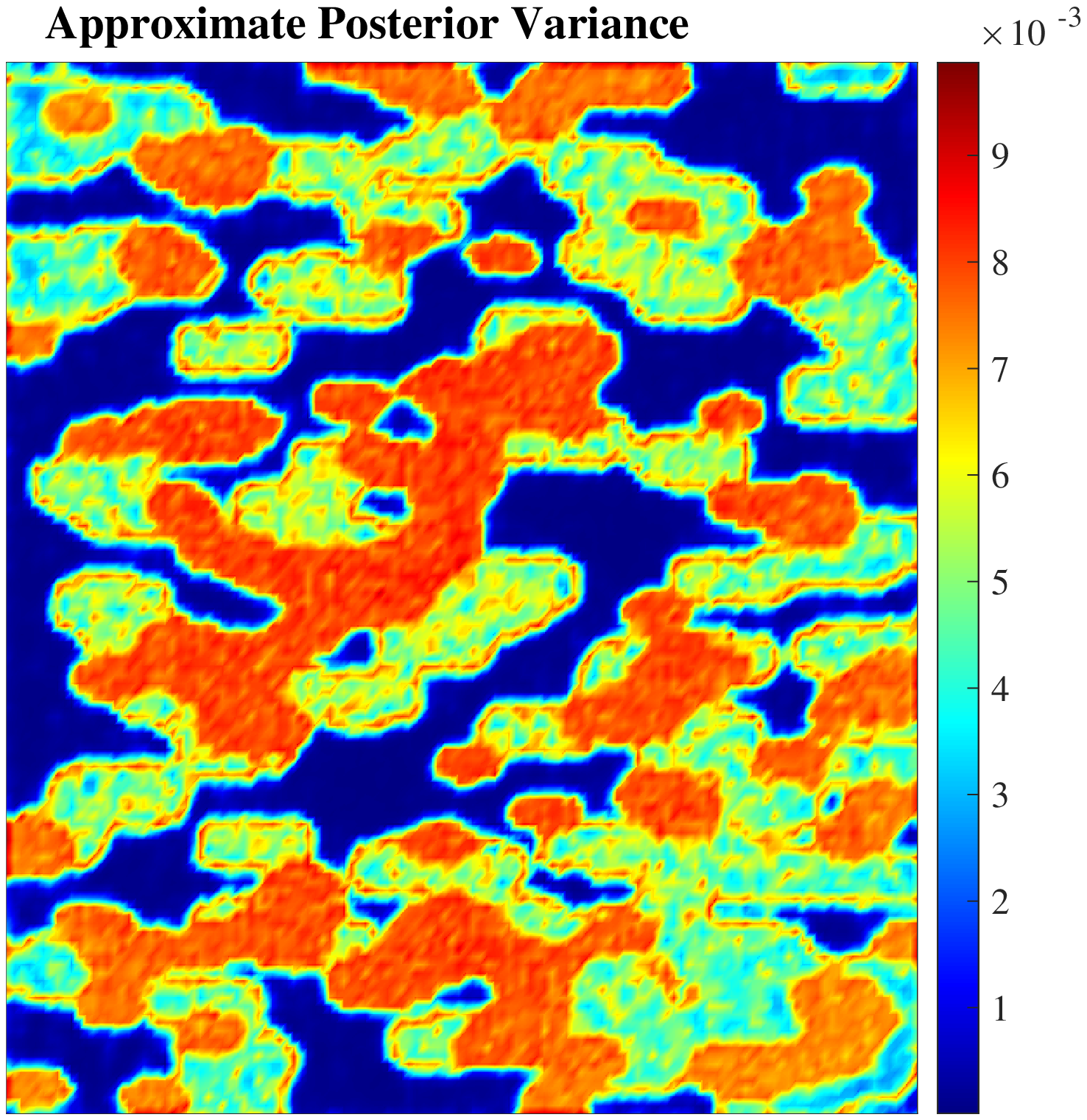}
    \caption{These results correspond to the `Three Phases' test problem: (left) Samples drawn from the approximate posterior distribution $\widehat\pi_{\vm}(\vm|\vd)$, and (right) approximate posterior variance computed using $N=1000$ Monte Carlo samples. }
    \label{fig:threephasesuq}
\end{figure}

 A summary of the relative error and the solver statistics for both test problems is provided in Table~\ref{tab:multiple_sum}. For the `Three Phases' test problems, we provide the true image and the MAP estimate in Figure~\ref{fig:threephasesmap}. In the bottom row of that same figure, we also provide the level set functions $\vphi_1$ and $\vphi_2$. In Figure~\ref{fig:threephasesuq}, we plot samples from the approximate posterior distribution $\widehat\pi_{\vm}(\vm|\vd)$, and on the right panel, we plot the approximate posterior variance computed using the LanczosMC approach (Method 1). We use $N=1000$ Monte Carlo samples and the size of the Lanczos basis was $k=200$. To solve linear systems with $\mh_{\rm GN}$, we used preconditioned GMRES (with no restart) instead of CG, even though the matrix is symmetric positive definite. This is because CG did not converge with this preconditioner in a reasonable number of iterations. The number of GMRES iterations, on average, were $273.3$. For generating the samples and for the LanczosMC approach, we used the Cholesky factorization of the $\mgamma_\text{prior}^{-1}$ as the preconditioner. The Lanczos solver for generating the samples took $409$ iterations, on average.

\begin{figure}[!ht]
    \centering
    \includegraphics[scale=0.4]{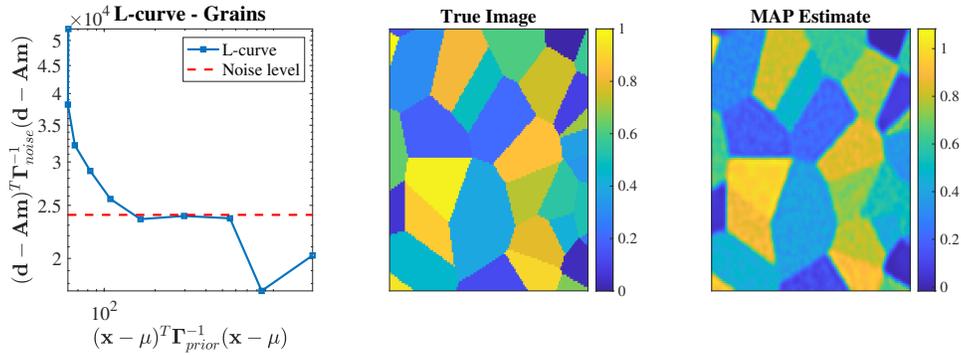}
    \caption{The image on the left displays the L-curve corresponding to noise levels $2\%$. The center image represents the true image and the image on the right displays the MAP estimate. The relative errors and a summary of the solver performance is given in Table~\ref{tab:multiple_sum}. }
    \label{fig:grainslcurve}
\end{figure}

\begin{figure}[!ht]
    \centering
    \includegraphics[scale=0.25]{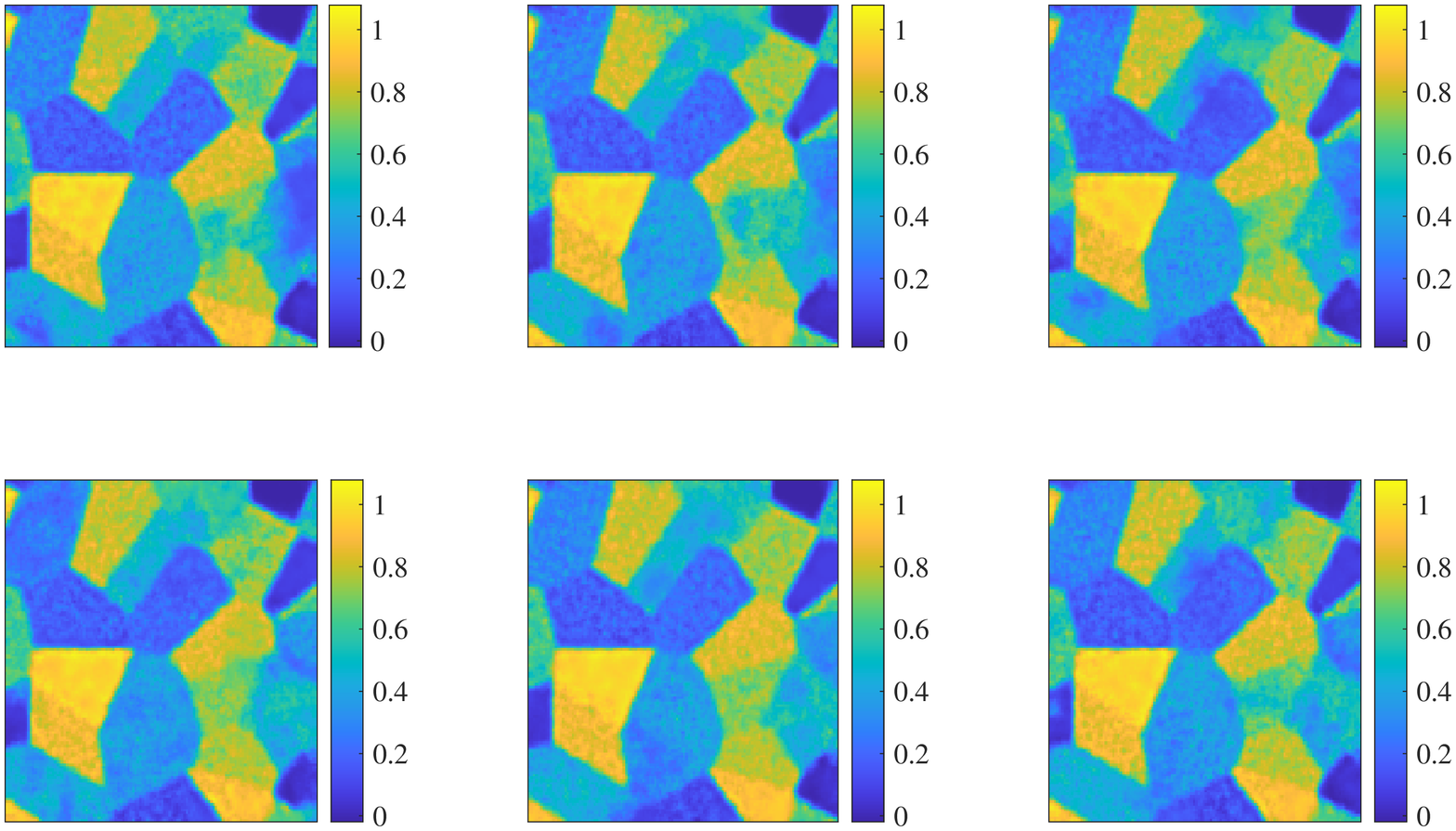}
    \includegraphics[scale=0.25]{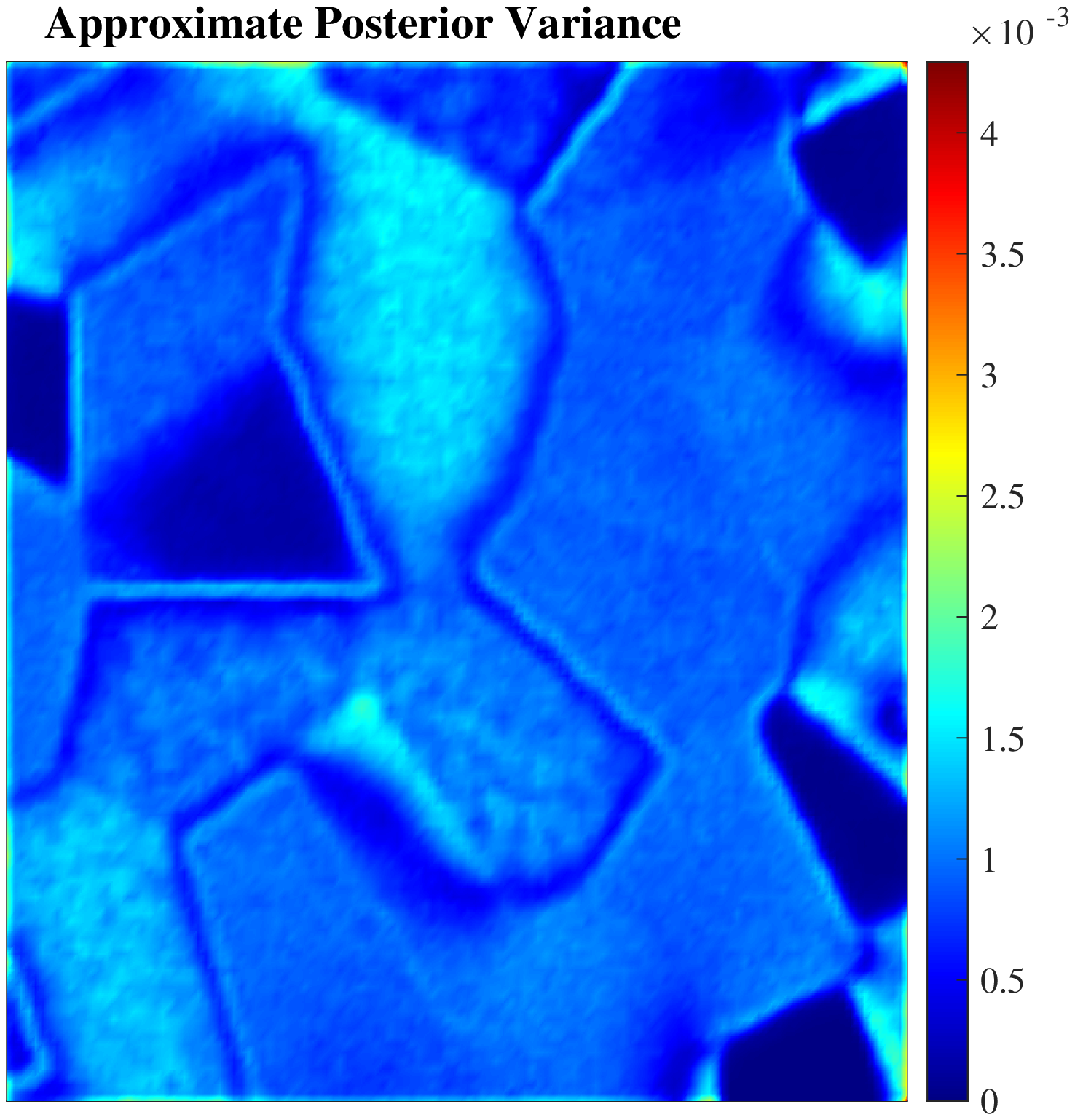}
    \caption{These results correspond to the `Grains' test problem corresponding to $2\%$ noise level: (top) Samples drawn from the Laplace's approximation to the posterior distribution $\widehat\pi_{\vm}(\vm|\vd)$, and (bottom) approximate posterior variance.} 
    \label{fig:grainsuq}
\end{figure}

 For the `Grains' test problems, we provide the true image and the MAP estimate in Figure~\ref{fig:grainslcurve}; we also plot the results from the L-curve for the grains example.  In Figure~\ref{fig:grainsuq}, we plot samples from the approximate posterior distribution $\widehat\pi_{\vm}(\vm|\vd)$, and on the right panel, we plot the approximate posterior variance computed using $N=1000$ Monte Carlo samples and the size of the Lanczos basis was $k=200$; the number of preconditioned GMRES iteration for solving systems with $\mh_{\rm GN}$ was $182.3$.  As before, for generating the samples, we used the Cholesky factorization of $\mgamma_\text{prior}^{-1}$ as the preconditioner and the Lanczos solver took $235.3$ iterations, on average.

\begin{table}[!ht] \centering
\begin{tabular}{c|c|c|c|c|c}
Image & $\lambda^2_{\vPhi} $ &	$ \text{relative error} $ & $\nls $ & \# GN & \# CG	\\
	\hline
Three Phases & $10^2/2^6$ &$ 13.9\%$ & $2$ & $8$ & $38$  \\
\hline
Grains & $10^3/2^5$ &$ 9.78\%$ & $3$ & $18$   & $1847$ 
\end{tabular}
\caption{Summary of the errors and the solver performance corresponding to `Three Phases' and `Grains' test problems. }
\label{tab:multiple_sum}
\end{table}

\subsubsection{Limited Angle Tomography} 
\begin{figure}[!ht]
    \centering
    \includegraphics[scale=0.5]{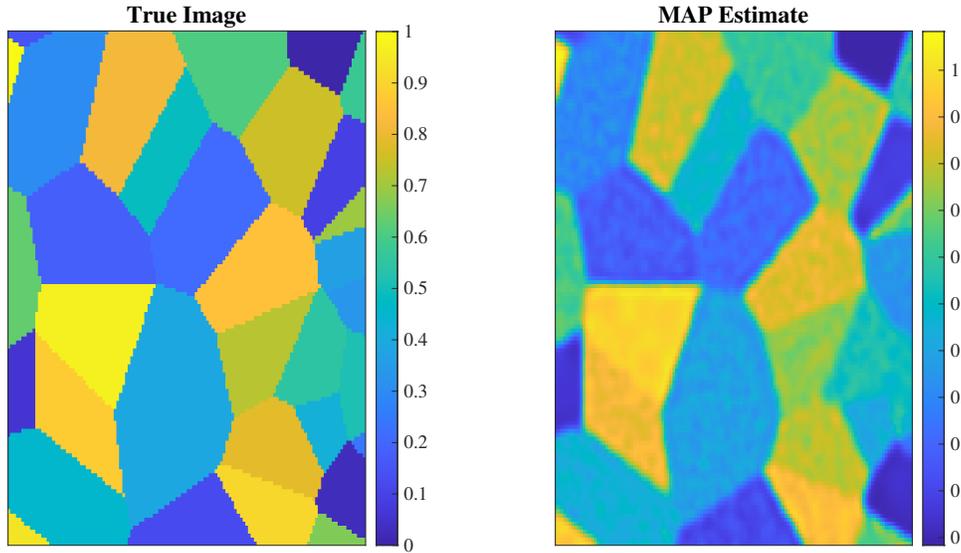}
    \caption{These images correspond to the `Grains' limited angle test problem: the true image (left) and the MAP estimate (right).}
    \label{fig:grainsmap_limited}
\end{figure}

For the next test problem, we consider a limited angle PAT problem. Measurements for the original forward operator $\ma$ for PAT  consists of line integrals along concentric circles whose centers span the original image~\cite{gazzola2019ir}.  For the limited angle case, we cut the number of angles in half, while maintaining the same angular range. As a consequence the new forward operator is under-determined and $n_{\rm obs} =11584 < n$. The test image we use is `Grains'. We set the noise level at $2\%$, and since the image has multiple regions we set $\nls = 3$ which can handle up to $2^{\nls}=8$ regions.  As before we take $\mgammac = \lambda^2_{\vc}\mi$ with $\lambda^2_{\vc} = 10\regPhi^2$.  In Figure \ref{fig:grainsmap_limited} we plot the true image and the corresponding map estimate. 

\begin{figure}[!ht]
    \centering
    \includegraphics[scale=0.23]{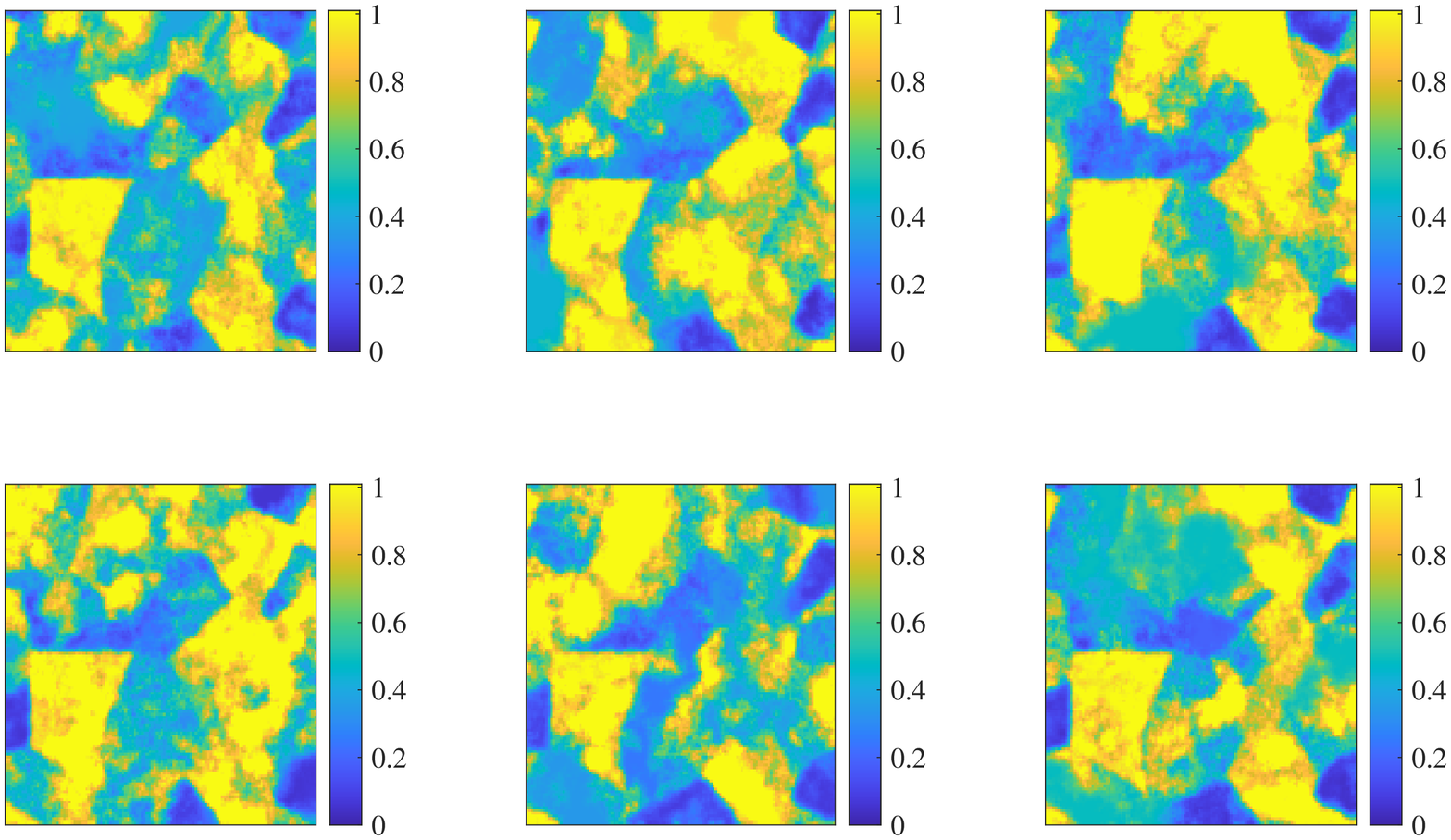}
    \includegraphics[scale=0.25]{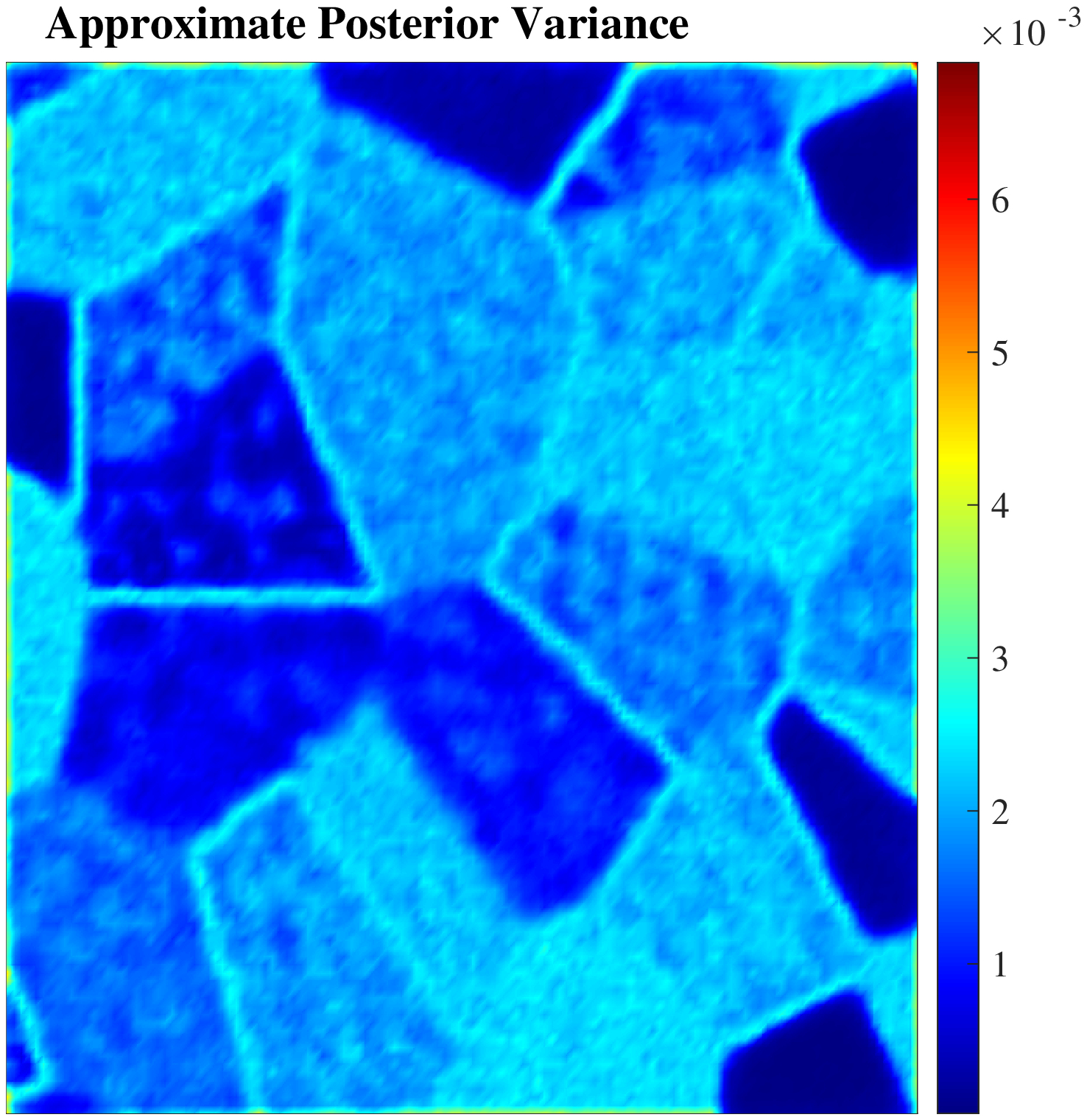}
    \caption{These results correspond to the `Grains' limited angle test problem with $2\%$ noise level and the limited angles: (left) Samples drawn from the approximate posterior distribution $\widehat\pi_{\vm}(\vm|\vd)$, and (right) approximate posterior variance computed using $N = 1000$ Monte Carlo samples. }
    \label{fig:grains_limiteduq}
\end{figure}

In Table~\ref{tab:grainsmap_limited} we summarize the error and solver statistics for the limited angle `Grains' example.  In Figure~\ref{fig:grains_limiteduq} we plot the approximate variance and samples from the posterior distribution for this example. The Lanczos solver for approximate samples took $203.8$ iterations, on average, and the number of GMRES iterations was $166.9$.

\begin{table}[!ht] \centering
\begin{tabular}{c|c|c|c|c|c}
noise & $\lambda^2_{\vPhi} $ &	$ \text{relative error} $ & $\nls $ & \# GN & \# CG	\\
\hline
$2 \%$ & $10^3/2^6$ &$ 11.2\%$ & $3$ & $6$  & $35$ 
\end{tabular}
\caption{Summary of the errors and the solver performance corresponding to limited angle `Grains' test problem.}
\label{tab:grainsmap_limited}
\end{table}

Looking at Tables~\ref{tab:multiple_sum} and~\ref{tab:grainsmap_limited} we notice that the number of CG iterations for the `Grains' example is much larger than for the `Three Phases example. One explanation for the `Grains' example having several more CG iterations could be a combination of its larger regularization parameter and the fact that it has a much smaller relative error than both 'Grains' limited angle and `Three phases.' Since the regularization parameter was larger it may have taken many more optimization iterations to converge.

\subsubsection{Comparing optimization solvers} To demonstrate the robustness of our inexact Gauss-Newton approach we will compare it to MATLAB's Trust Region algorithm and the LBFGS algorithm from the Poblano toolbox \cite{SAND2010-1422}. By default, the LBFGS solver has a limited memory parameter of $5$. The default stopping criterion for the LBFGS algorithm are: a maximum of 100 optimization iterations, a maximum of 100 function evaluations, a scaled gradient norm less than or equal to $10^{-5}$, and a relative function change less than or equal to $10^{-6}$. See \cite{SAND2010-1422} for more details. As discussed in Subsection \ref{ssec:gn} our solver (inexact Gauss-Newton) uses both the discrepancy principle and a relative reduction in the gradient as stopping criterion. For the gradient stopping criterion we require that the initial gradient reduces by three orders of magnitude. For all of the optimizers we will set the maximum number of optimization iterations to $50$. We set the maximum number of CG iterations per optimization iteration to $200$ and set the tolerance for the relative residuals to $10^{-6}$. 
Note that we provide the Gauss-Newton Hessian in lieu of the exact Hessian for the Trust Region solver. We will examine the performance of computing the map estimate on the `Grains' test problem. For this example we use the same noise and  $\regPhi^2$ as reported in Table~\ref{tab:multiple_sum}. To compare we report the number of objective function evaluations, optimization iterations, the total number of CG iterations and report the relative error. Note that CG iterations are not reported for the LBFGS solver because the search direction in the LBFGS algorithm is calculated using a Hessian inverse approximation formula~\cite{SAND2010-1422}. The results are given in Table~\ref{tab:optimizers}.

\begin{table}[H] \centering
\begin{tabular}{c|c|c|c|c}
    Solver & \#F evals & \#Opt iters & \#CG iters & Rel error\\ \hline
    LBFGS & 120 & 50$^{(*)}$ & N/A & 11.4 \% \\ \hline 
    Trust Region & 52 & 50$^{(*)}$ & 8200 & 9.68 \% \\ \hline
    Inexact GN  & 71 & 18 & 1847 & 9.78\% \\ 
     \end{tabular}
     \caption{Table comparing performance of different optimization solvers for the Grains example. The symbol $(*)$ indicates that the optimizer did not converge within the $50$ iteration limit. }
     \label{tab:optimizers}
\end{table}
Table \ref{tab:optimizers}  demonstrates the benefits of utilizing an inexact Gauss-Newton approach in favor of other popular methods such as LBFGS. LBFGS completed $120$ function evaluations as opposed to Trust Region and both inexact Gauss-Newton methods which were $52$ and $71$ respectively. The Trust Region method utilized less function evaluations than the inexact Gauss-Newton  methods but had almost four times the number of CG iterations. One explanation for this may be that since inexact Gauss-Newton only requires CG to be solved within a tolerance specified by the forcing sequence as in~\eqref{eq:forcingSequence}. On the other hand, inexact Gauss-Newton requires more function evaluations at the linesearch stage.

\subsection{X-ray tomography}
In this set of numerical experiments, we use the real-data from~\cite{hamalainen2015tomographic}. The dataset consists of the logarithm of the x-ray sinogram of a single two-dimensional slice of a walnut and the forward operator $\ma$ modeling the x-ray transform for this particular set.  We use all 120 projections that are provided in the dataset; therefore, the size of the matrix $\ma$ is $ 9840\times 6724$. From the ground truth image, while the image is not exactly piecewise constant, we use $\nls = 2$ which yields $2^{\nls} = 4$ regions. The parameter $\lambda_{\vPhi}^2$ is chosen using an L-curve type analysis; although we do not know the exact noise level we use the curvature of the L-curve to guide the selection of the regularization parameter. We take the noise covariance matrix to be $\mgamma_\mathrm{noise} = \sigma^2 \mi$ where $\sigma = 0.035$ and $\lambda_{\vPhi} = 10^3/2^5$. We also choose $\lambda_{\vc}^2 = 50 \lambda_{\vPhi}^2$.

\begin{figure}[!ht]
    \centering
    \includegraphics[scale=0.3]{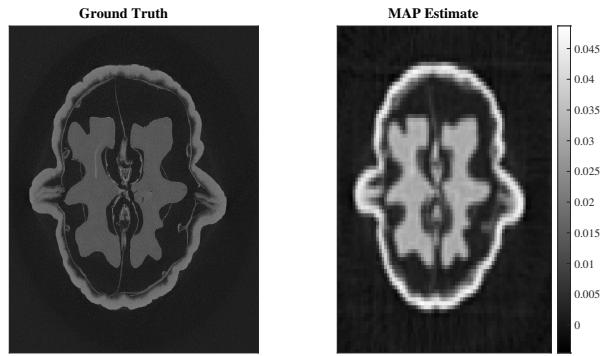}
    \caption{(left) Ground truth obtained using the filtered backprojection technique. The image is provided at a much higher resolution and is meant for visual comparison only. (right) MAP estimate obtained using $\nls =2$ level set functions. }
    \label{fig:walnutmap}
\end{figure}

\begin{figure}[!ht]
    \centering
    \includegraphics[scale=0.25]{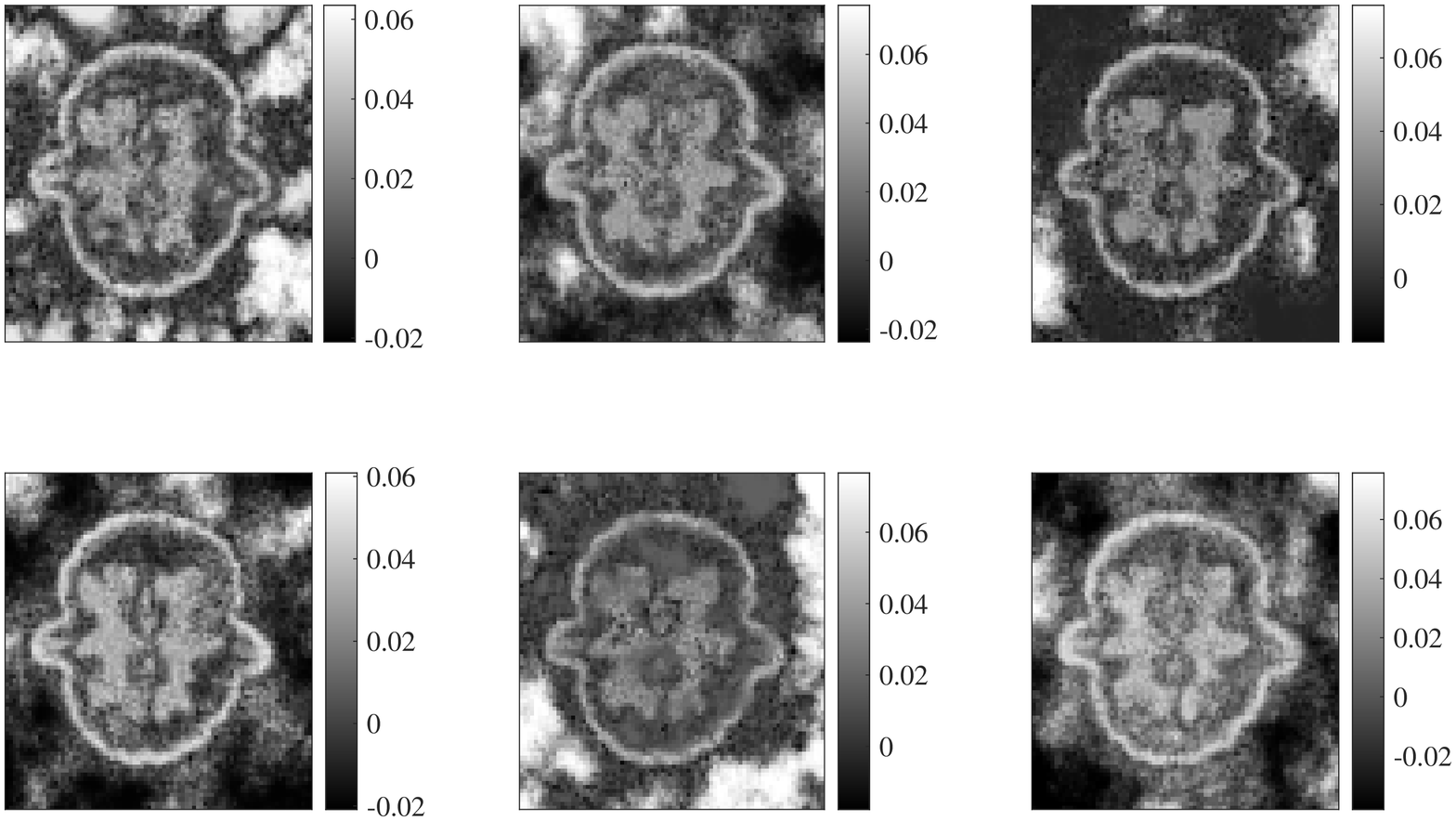}
    \includegraphics[scale=0.25]{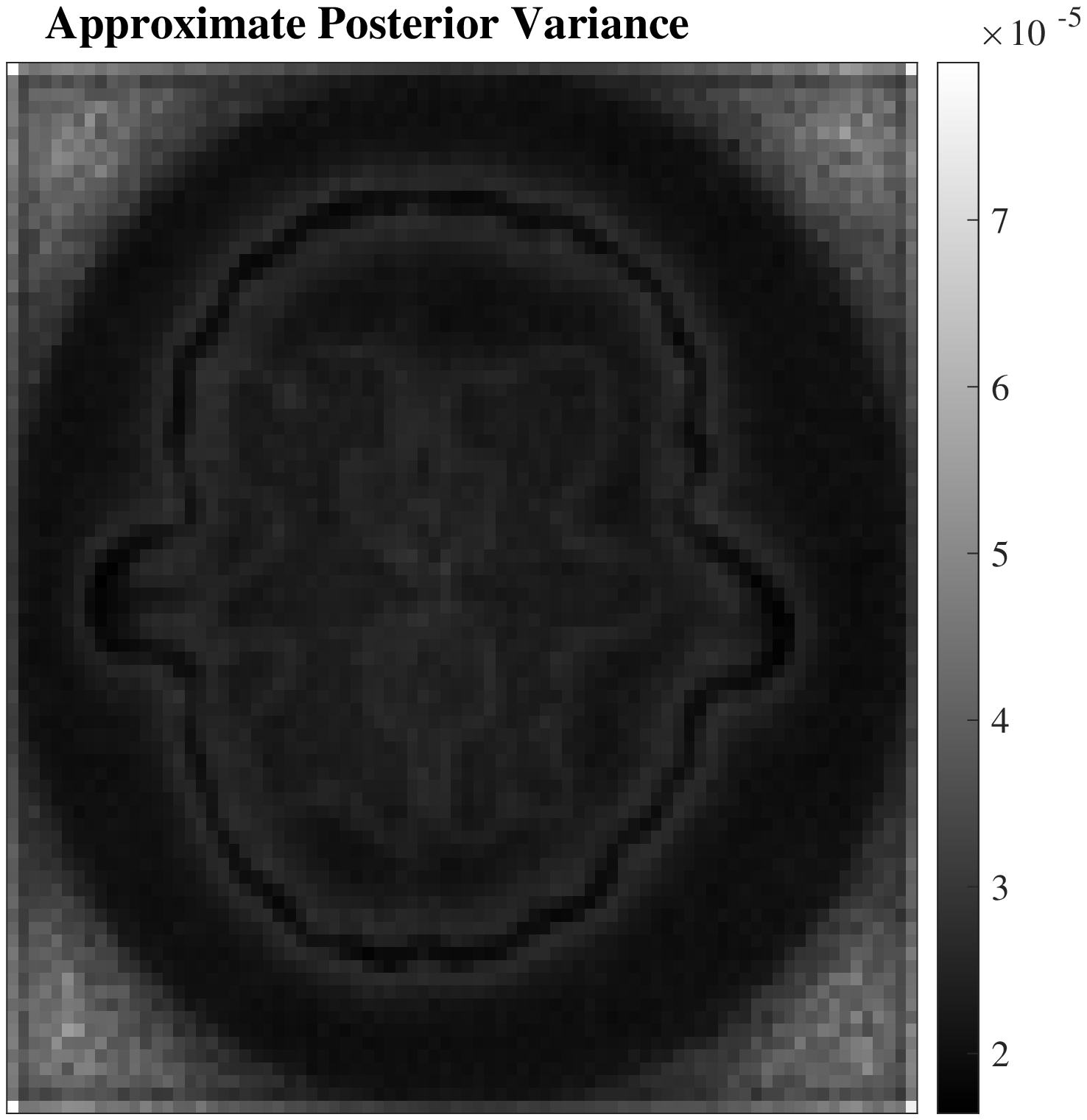}
    \caption{These results correspond to the `Walnut' test problem: (left) Samples drawn from the approximate posterior distribution $\widehat\pi_{\vm}(\vm|\vd)$, and (right) approximate posterior variance. }
    \label{fig:walnutuq}
\end{figure}

To compute the MAP estimate, we use the inexact Gauss-Newton solver which converged in $18$ iterations and required $683$ CG iterations overall. Since prior knowledge of the magnitude of the noise is not available, the optimization routine uses a stopping criterion based on the reduction of the relative magnitude by four orders of magnitude. The MAP estimate has been displayed alongside the ground truth image obtained using a filtered backprojection technique; note that this image is at a much higher resolution than what is being used to solve the inverse problem. From a visual perspective, it is clear that the MAP estimate captures the main outline of the walnut shape and we do not provide a quantitative comparison for the reasons described above. In addition to the MAP estimate, we provide samples from the approximate posterior $\widehat\pi_{\vm}(\vm|\vd)$ in the left panel of Figure~\ref{fig:walnutuq}. Finally, we provide the approximate posterior variance, computed using $N=1000$ Monte Carlo samples and Lanczos basis $k=100$, in the right panel of the same figure. The number of GMRES iterations, on average, was $195.8$. The Lanczos solver for the approximate samples took $248.7$ iterations, on average. The image shows a higher variance around the edges of the walnut and at the ends of the images.  

\subsection{Hydraulic Tomography}
\begin{figure}[!ht]
    \centering
    \includegraphics[scale=0.4]{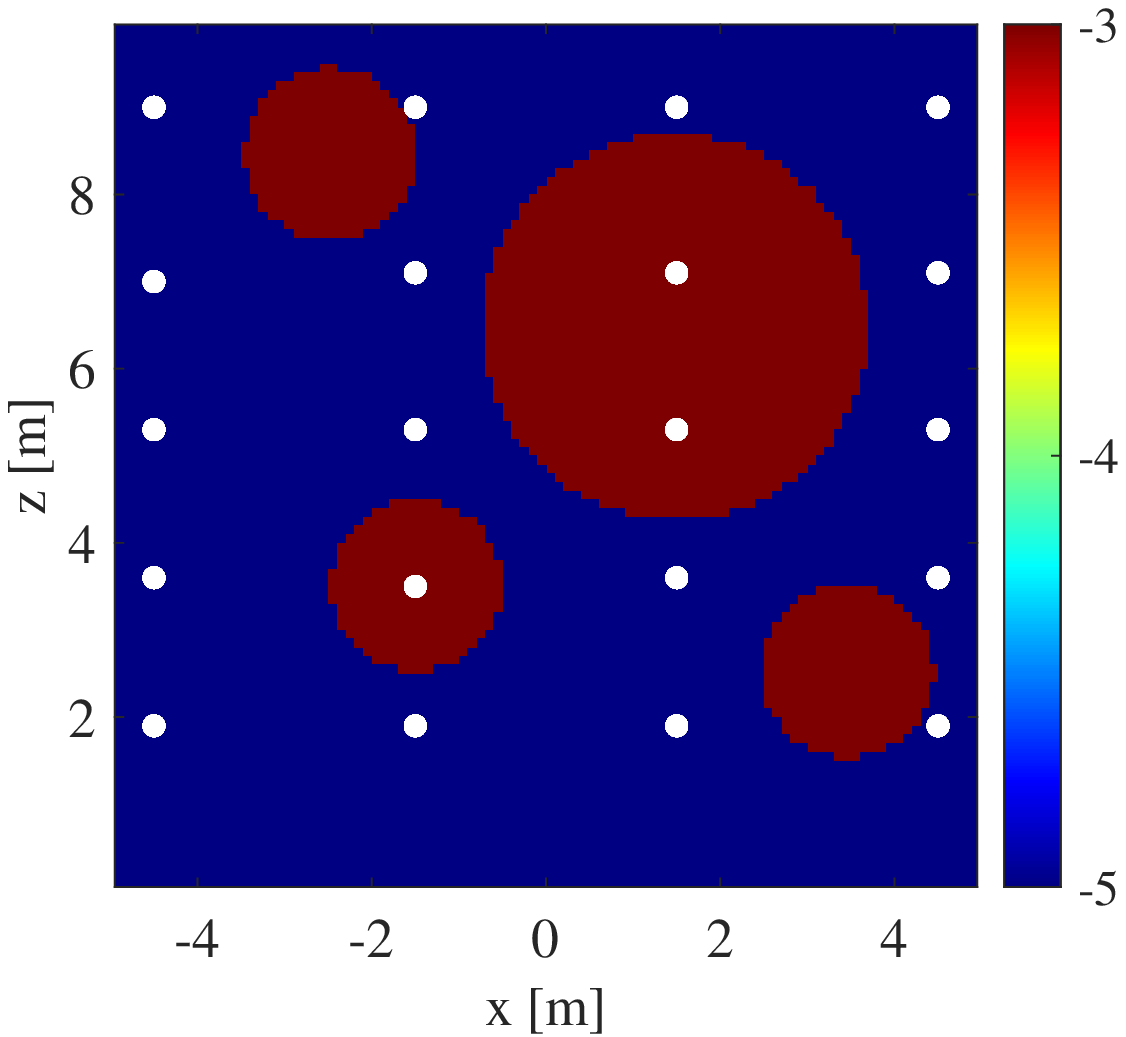}
    \includegraphics[scale=0.3]{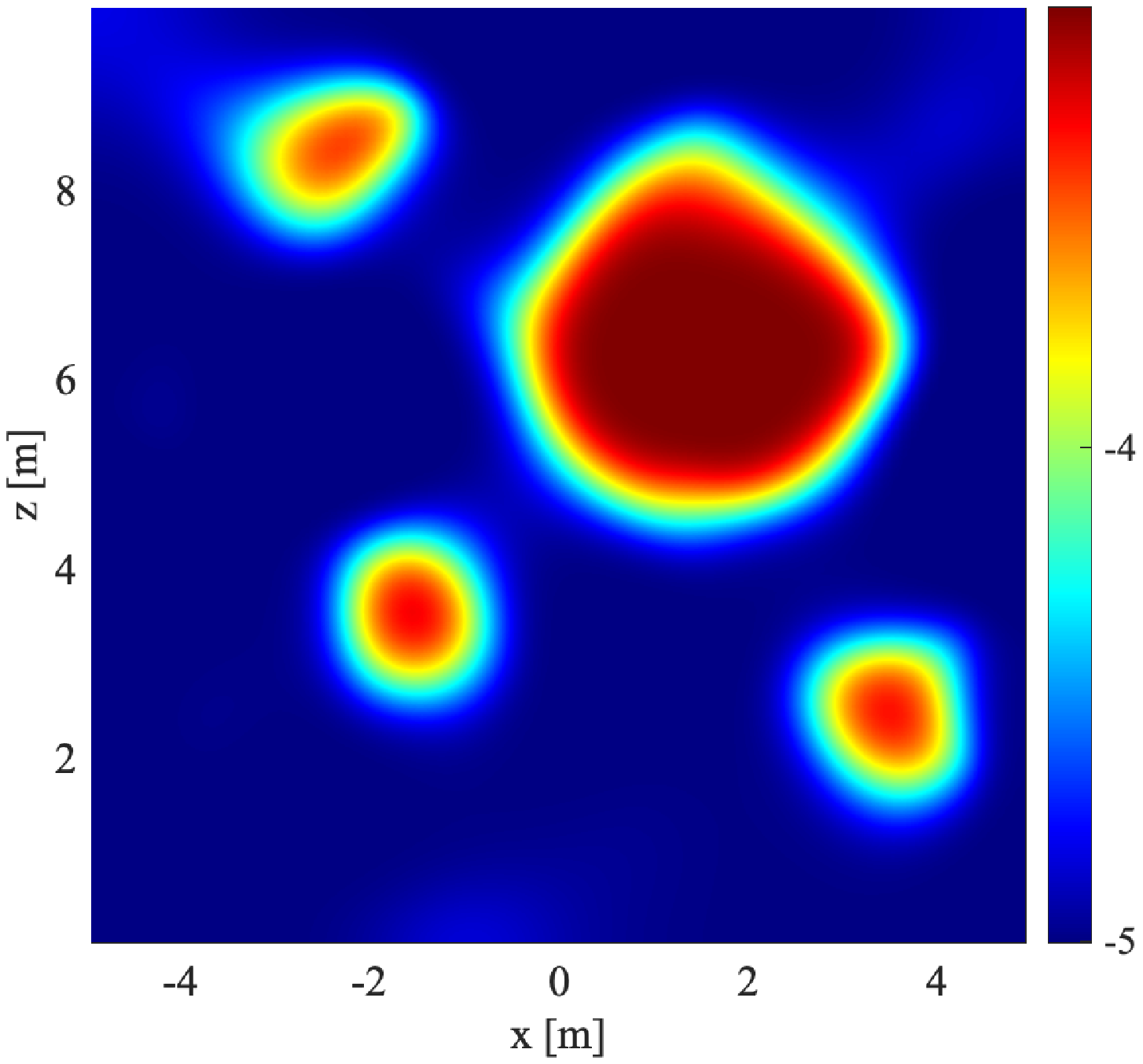}
    \caption{Ground truth log hydraulic conductivity field ($\log_{10} K$) for a hydraulic tomography example consisting of gravel ($\log_{10} K = -5$ m/s) and sand ($\log_{10} K = -3$ m/s). White circles are the pumping/observation well locations.}
    \label{fig:HT}
\end{figure}

\begin{figure}[!ht]
    \centering
    \includegraphics[scale=0.3]{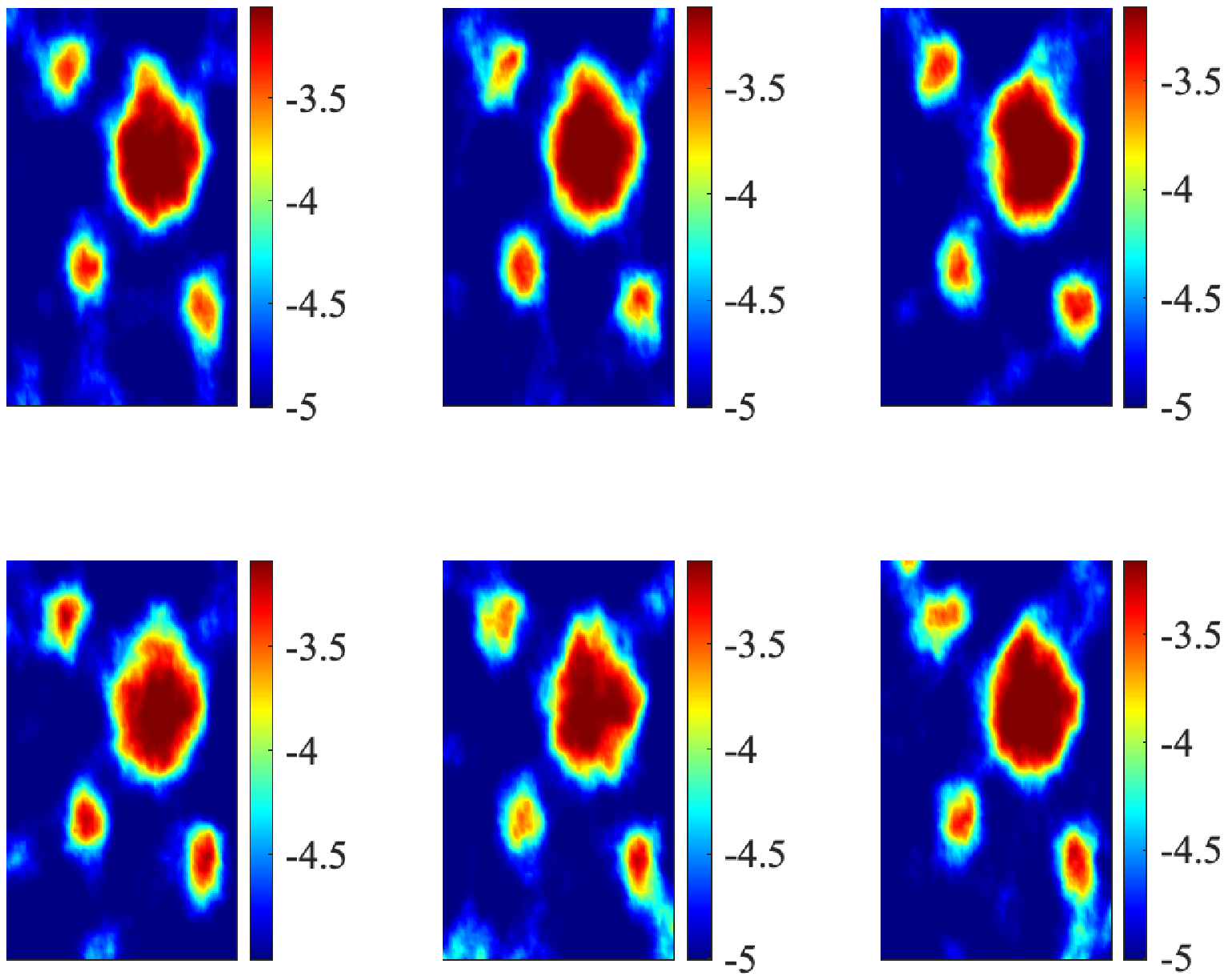}
    \includegraphics[scale=0.3]{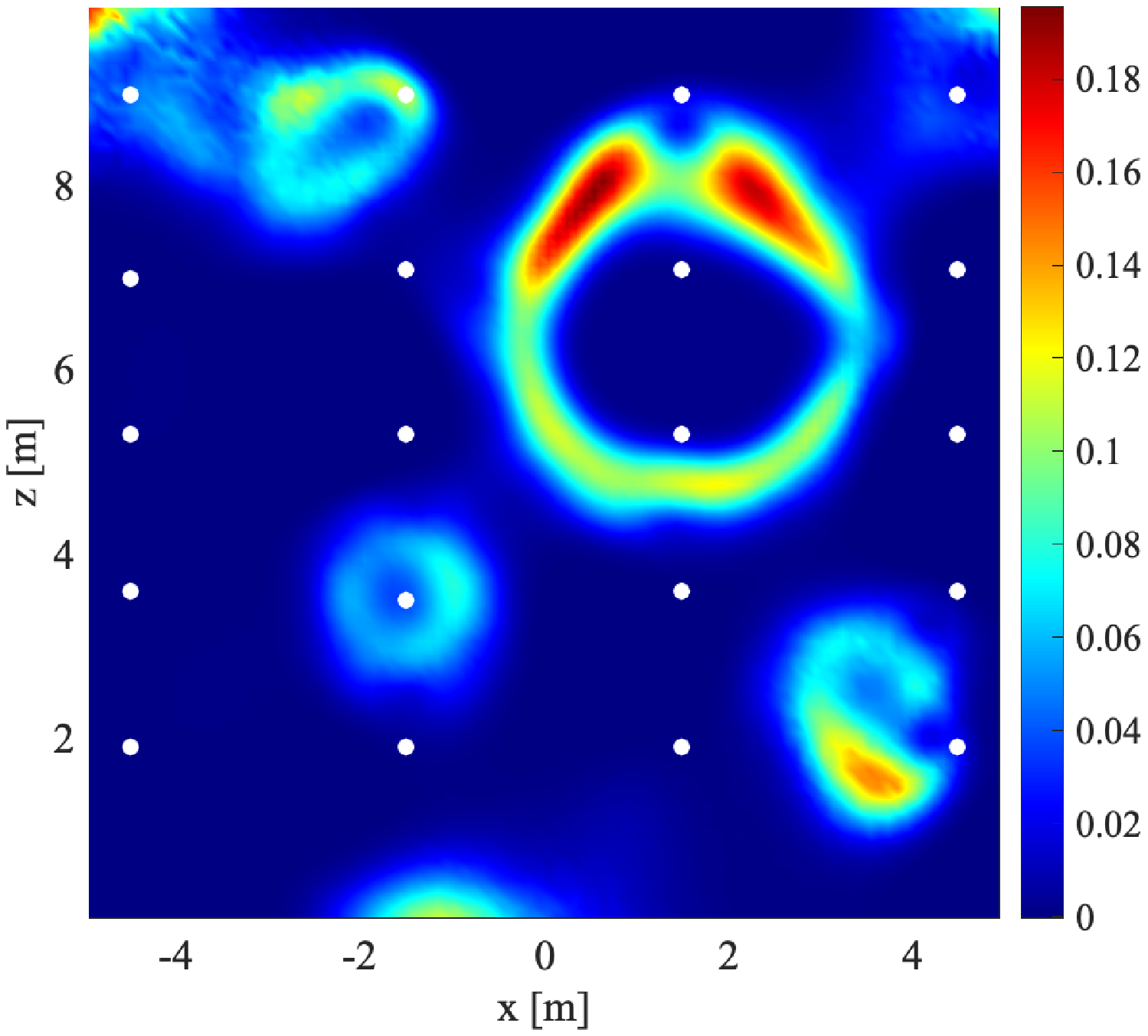}
    \caption{These results correspond to the hydraulic tomography test problem: (left) Samples from the approximate posterior distribution $\widehat\pi_{\vm}(\vm|\vd)$ and (right) approximate posterior variance. White circles are the pumping/observation well locations.}
    \label{fig:HT_UQ}
\end{figure}

In this example, we present a nonlinear hydraulic tomography problem adapted from~\cite{lee2013}. A synthetic vertical confined aquifer consisting of two geologic facies is considered and the subsurface flow system is governed by
\begin{equation}\label{eq:gw}
    -\nabla \cdot \left(K \nabla h_i \right) = q_i \delta\left(\vs_i \right) \qquad 1 \leq i \leq n_s 
\end{equation}
where $K$ is the hydraulic conductivity [$m/s$], $h$ is the hydraulic head [$m$], $q_i$ is the 2-D pumping rate for a pumping test at a well location $\vs_i$ (marked with a circle in Figure~\ref{fig:HT}) for $1 \leq i \leq n_s$, $n_s$  is the number of sources, and $\delta(\vx)$ is the Dirac delta function. A constant-head boundary condition on the left, right, and upper boundaries as well as no-flux condition on the bottom is imposed. The measurements are steady state head changes due to sequential pump tests at each well location. For each test water is extracted from one of the $20$ well locations at the rate of $2$ L/s and the change in the head is measured at all other 19 locations resulting in $380$ drawdown measurements ($20$ pump tests $\times$ $19$ measurements per test). To ensure positive hydraulic conductivity values during the inversion, log-transformed hydraulic conductivity $\log K$ is estimated and 1\% Gaussian noise ($\approx$ std(error) = $0.12$ m) was added to the measurements. The domain was discretized into $100$ by $100$ grids so that $n=10^4$. The governing equation in~\eqref{eq:gw} is solved by a finite volume method.    
To compute the MAP estimate, the inexact Gauss-Newton solver is used and it converged in $17$ iterations with $1006$ CG iterations in total. Because of log-transformation, the prior covariance for $\vPhi$ is  $\mgammaPhi^{-1} = \mi_{\nls} \otimes  \regPhi^2\left(\alpha\ml + \gamma\mi \right)^2$ and for $\vc$ as a Gaussian distribution with mean $\log_{10} K = -4$ and covariance matrix $\mgammac=\lambda_{\vc}^2\mi$ and $\lambda_{\vc}^2 =  \lambda_{\vPhi}^2$. The MAP estimate is displayed along with the true  $\log_{10}K$ field in Figure~\ref{fig:HT}. The estimate identifies the main circle structures with different values of $K$. 
Similar to the previous examples we plot samples from the approximate posterior distribution as well as the approximate posterior variance (using LanczosMC with the same settings as before) in Figure~\ref{fig:HT_UQ}. The Lanczos solver took $44.2$ iterations, on average. To compute the approximate posterior variance, we used preconditioned CG instead of GMRES in contrast to the previous applications; the total number of iterations, on average, was $207.9$.

In Figure~\ref{fig:HT_UQ} we see that the magnitude of the  uncertainty is high at the interfaces of the regions especially in the areas far away from the well locations as expected. The high uncertainty regions are well matched with the interfaces of the geologic structures. In the interior of largest high permeable circle there is a decrease in uncertainty. The largest circle contains more measurements than the other ones which would imply that the uncertainty would be smaller there. 
\section{Conclusions}\label{sec:conc}
In this paper we have developed a Bayesian level set approach to inverse problems where the unknown is a piecewise constant field. This was accomplished by assigning prior distributions to level set functions which define different piecewise constant regions and the constants that define their magnitudes. Assuming an additive  Gaussian noise model and Gaussian prior we derived the (potentially) non-Gaussian posterior distribution and developed an inexact Gauss-Newton approach to compute the MAP estimate that leveraged efficient adjoint-based methods for computing derivatives and sparsity in the Jacobian,  of the piecewise constant model. We employed Laplace's approximation to approximate the posterior and developed matrix-free methods for sampling from the posterior and estimating the posterior variance. In addition we developed a matrix-free estimator for the diagonals of the inverse of a matrix. Our algorithms were performed on computational examples in photoacoustic, hydraulic, and X-ray tomography. In each example we obtained an accurate reconstruction and gained qualitative information about the uncertainty in the reconstruction using samples from the approximate posterior and the approximate posterior variance. It is important to emphasize that the approximated uncertainty estimates are to be interpreted qualitatively.  Although we have proposed  efficient and matrix-free methods to Bayesian inverse problems with piecewise constant unknowns there are a number of possible future directions for our work.

An immediate extension of this work would be to implement a Gaussian prior with the Mat\'{e}rn covariance. The Mat\'{e}rn covariance is a stationary covariance kernel that allows for control of the smoothness and variance of prior samples \cite{matern}. As a result, the Mat\'{e}rn covariance provides more flexibility than the inverse elliptical differential operator we have implemented in this work. This flexibility comes at the cost of sparsity in the prior precision matrix $\mgamma^{-1}_\text{prior}$ which we have leveraged in our algorithms. The new prior covariance will be a large dense matrix and explicitly forming or factorizing the precision matrix will be computationally infeasible and introduces new challenges in our approaches for MAP estimation and uncertainty quantification. To tackle these challenges we would build on the methods in \cite{genhybr}. A topic that we have not addressed in this paper is the validation of linearized Bayesian inference. One way to accomplish this would be to compare samples generated from the method in~\ref{ssec:sample} to samples generated from a Gauss-Newton Hessian based MCMC sampler \cite{petra2014computational}.

\section{Acknowledgements} This work was supported, in part, by the National Science Foundation through the award DMS-1845406 (W.R. and A.K.S.) and DGE-1746939 (W.R.).

\section{References}
\bibliography{refs}

\begin{thebibliography}{10}

\bibitem{aghasi2011parametric}
A.~Aghasi, M.~Kilmer, and E.~L. Miller.
\newblock Parametric level set methods for inverse problems.
\newblock {\em SIAM Journal on Imaging Sciences}, 4(2):618--650, 2011.

\bibitem{bardsley2018computational}
J.~M. Bardsley.
\newblock {\em Computational uncertainty quantification for inverse problems},
  volume~19 of {\em Computational Science \& Engineering}.
\newblock Society for Industrial and Applied Mathematics (SIAM), Philadelphia,
  PA, 2018.

\bibitem{bekas2007estimator}
C.~Bekas, E.~Kokiopoulou, and Y.~Saad.
\newblock An estimator for the diagonal of a matrix.
\newblock {\em Applied numerical mathematics}, 57(11-12):1214--1229, 2007.

\bibitem{bui2015scalable}
T.~Bui-Thanh and O.~Ghattas.
\newblock A scalable algorithm for {MAP} estimators in {B}ayesian inverse
  problems with {B}esov priors.
\newblock {\em Inverse Problems \& Imaging}, 9(1):27, 2015.

\bibitem{buithanh2013computational}
T.~Bui-Thanh, O.~Ghattas, J.~Martin, and G.~Stadler.
\newblock A computational framework for infinite-dimensional {B}ayesian inverse
  problems part {I}: {T}he linearized case, with application to global seismic
  inversion.
\newblock {\em SIAM Journal on Scientific Computing}, 35(6):A2494--A2523, 2013.

\bibitem{burger2005survey}
M.~Burger and S.~J. Osher.
\newblock A survey on level set methods for inverse problems and optimal
  design.
\newblock {\em European journal of applied mathematics}, 16(2):263--301, 2005.

\bibitem{blip}
M.~Cardiff and P.~Kitanidis.
\newblock Bayesian inversion for facies detection: An extensible level set
  framework.
\newblock {\em Water Resources Research}, 45, 2009.

\bibitem{chantas2008}
G.~Chantas, N.~Galatsanos, A.~Likas, and M.~Saunders.
\newblock Variational {B}ayesian image restoration based on a product of
  t-distributions image prior.
\newblock {\em IEEE Transactions on Image Processing}, 17(10):1795--1805, 2008.

\bibitem{chow2014preconditioned}
E.~Chow and Y.~Saad.
\newblock Preconditioned {K}rylov subspace methods for sampling multivariate
  {G}aussian distributions.
\newblock {\em SIAM Journal on Scientific Computing}, 36(2):A588--A608, 2014.

\bibitem{suitesparse}
T.~A. Davis and Y.~Hu.
\newblock The university of {F}lorida sparse matrix collection.
\newblock {\em ACM Transactions on Mathematical Software}, 38(1):1--25, 2011.

\bibitem{dorn2006level}
O.~Dorn and D.~Lesselier.
\newblock Level set methods for inverse scattering.
\newblock {\em Inverse Problems}, 22(4):R67, 2006.

\bibitem{SAND2010-1422}
D.~M. Dunlavy, T.~G. Kolda, and E.~Acar.
\newblock Poblano v1.0: A matlab toolbox for gradient-based optimization.
\newblock Number SAND2010-1422, March 2010.

\bibitem{dunlop2017hierarchical}
M.~M. Dunlop, M.~A. Iglesias, and A.~M. Stuart.
\newblock Hierarchical {B}ayesian level set inversion.
\newblock {\em Statistics and Computing}, 27(6):1555--1584, 2017.

\bibitem{dunlop2016map}
M.~M. Dunlop and A.~M. Stuart.
\newblock {MAP} estimators for piecewise continuous inversion.
\newblock {\em Inverse Problems}, 32(10):105003, 2016.

\bibitem{gazzola2019ir}
S.~Gazzola, P.~C. Hansen, and J.~G. Nagy.
\newblock {IR} tools: a {MATLAB} package of iterative regularization methods
  and large-scale test problems.
\newblock {\em Numerical Algorithms}, 81(3):773--811, 2019.

\bibitem{matern}
M.~G. Genton.
\newblock Classes of kernels for machine learning: A statistics perspective.
\newblock {\em J. Mach. Learn. Res.}, 2:299–312, mar 2002.

\bibitem{haber2000optimization}
E.~Haber, U.~M. Ascher, and D.~Oldenburg.
\newblock On optimization techniques for solving nonlinear inverse problems.
\newblock {\em Inverse problems}, 16(5):1263, 2000.

\bibitem{hamalainen2015tomographic}
K.~H{\"a}m{\"a}l{\"a}inen, L.~Harhanen, A.~Kallonen, A.~Kujanp{\"a}{\"a},
  E.~Niemi, and S.~Siltanen.
\newblock Tomographic x-ray data of a walnut.
\newblock {\em arXiv preprint arXiv:1502.04064}, 2015.

\bibitem{hansen2010discrete}
P.~C. Hansen.
\newblock {\em Discrete inverse problems}, volume~7 of {\em Fundamentals of
  Algorithms}.
\newblock Society for Industrial and Applied Mathematics (SIAM), Philadelphia,
  PA, 2010.
\newblock Insight and algorithms.

\bibitem{iglesias2016}
M.~A. Iglesias, Y.~Lu, and A.~M. Stuart.
\newblock A {B}ayesian level set method for geometric inverse problems.
\newblock {\em Interfaces and Free Boundaries}, 18, 2016.

\bibitem{kaipio2005statistical}
J.~Kaipio and E.~Somersalo.
\newblock {\em Statistical and computational inverse problems}, volume 160 of
  {\em Applied Mathematical Sciences}.
\newblock Springer-Verlag, New York, 2005.

\bibitem{siltanen2004}
M.~Lassas and S.~Siltanen.
\newblock Can one use total variation prior for edge preserving {B}ayesian
  inversion?
\newblock {\em Inverse Problems}, 20(5):1537--1563, 2004.

\bibitem{lassi2014whittle}
S.~L. Lassi~Roininen, Janne M. J.~Huttunen.
\newblock {W}hittle-{M}at\'ern priors for {B}ayesian statistical inversion with
  applications in electrical impedance tomography.
\newblock {\em Inverse Problems \& Imaging}, 8(2):561--586, 2014.

\bibitem{lee2013}
J.~Lee and P.~Kitandis.
\newblock Bayesian inversion with total variation prior for discrete geologic
  structure identification.
\newblock {\em Water Resources Research}, 49:7658--7669, 2013.

\bibitem{lehmann2006theory}
E.~L. Lehmann and G.~Casella.
\newblock {\em Theory of point estimation}.
\newblock Springer Science \& Business Media, 2006.

\bibitem{lindgren2011explicit}
F.~Lindgren, H.~Rue, and J.~Lindstr{\"o}m.
\newblock An explicit link between {G}aussian fields and {G}aussian {M}arkov
  random fields: the stochastic partial differential equation approach.
\newblock {\em Journal of the Royal Statistical Society: Series B (Statistical
  Methodology)}, 73(4):423--498, 2011.

\bibitem{MarkkanenRoininenHuttunenLasanen+2019+225+240}
M.~Markkanen, L.~Roininen, J.~M. Huttunen, and S.~Lasanen.
\newblock Cauchy difference priors for edge-preserving {B}ayesian inversion.
\newblock {\em Journal of Inverse and Ill-posed Problems}, 27(2):225--240,
  2019.

\bibitem{meyer2021hutch++}
R.~A. Meyer, C.~Musco, C.~Musco, and D.~P. Woodruff.
\newblock Hutch++: Optimal stochastic trace estimation.
\newblock In {\em Symposium on Simplicity in Algorithms (SOSA)}, pages
  142--155. SIAM, 2021.

\bibitem{nocedal1999}
J.~Nocedal and S.~J. Wright.
\newblock {\em Numerical Optimization}.
\newblock Springer, 1 edition, 1999.

\bibitem{petra2014computational}
N.~Petra, J.~Martin, G.~Stadler, and O.~Ghattas.
\newblock A computational framework for infinite-dimensional {B}ayesian inverse
  problems, part ii: {S}tochastic {N}ewton {MCMC} with application to ice sheet
  flow inverse problems.
\newblock {\em SIAM Journal on Scientific Computing}, 36(4):A1525--A1555, 2014.

\bibitem{rasmussen2003gaussian}
C.~E. Rasmussen.
\newblock Gaussian processes in machine learning.
\newblock In {\em Summer school on machine learning}, pages 63--71. Springer,
  2003.

\bibitem{GMRFbook}
H.~Rue and L.~Held.
\newblock {\em Gaussian {M}arkov Random Fields: {T}heory and Applications},
  volume 104 of {\em Monographs on Statistics and Applied Probability}.
\newblock Chapman \& Hall, London, 2005.

\bibitem{genhybr}
A.~Saibaba, J.~Chung, and K.~Petroske.
\newblock Efficient {K}rylov subspace methods for uncertainty quantification in
  large {B}ayesian linear inverse problems.
\newblock {\em Numerical Linear Algebra with Applications}, 27, 08 2020.

\bibitem{stuart2010inverse}
A.~M. Stuart.
\newblock Inverse problems: a {B}ayesian perspective.
\newblock {\em Acta Numerica}, 19:451--559, 2010.

\bibitem{van2010multiple}
K.~van~den Doel, U.~M. Ascher, and A.~Leitao.
\newblock Multiple level sets for piecewise constant surface reconstruction in
  highly ill-posed problems.
\newblock {\em Journal of Scientific Computing}, 43(1):44--66, 2010.

\bibitem{villa2021hippylib}
U.~Villa, N.~Petra, and O.~Ghattas.
\newblock {hIPPYlib: An Extensible Software Framework for Large-Scale Inverse
  Problems Governed by PDEs: Part I: Deterministic Inversion and Linearized
  Bayesian Inference}.
\newblock {\em ACM Transactions on Mathematical Software (TOMS)}, 47(2):1--34,
  2021.

\bibitem{vogel2002computational}
C.~R. Vogel.
\newblock {\em Computational methods for inverse problems}, volume~23 of {\em
  Frontiers in Applied Mathematics}.
\newblock Society for Industrial and Applied Mathematics (SIAM), Philadelphia,
  PA, 2002.
\newblock With a foreword by H. T. Banks.

\bibitem{wu2016estimating}
L.~Wu, J.~Laeuchli, V.~Kalantzis, A.~Stathopoulos, and E.~Gallopoulos.
\newblock Estimating the trace of the matrix inverse by interpolating from the
  diagonal of an approximate inverse.
\newblock {\em Journal of Computational Physics}, 326:828--844, 2016.

\end{thebibliography}
\bibliographystyle{abbrv}

\appendix
\section{Estimating the Diagonal of $\ma^{-1}$}\label{sec:diaginv}

In this appendix, we give the details of the new estimator (LanczosMC) for the diagonals of the inverse of  a positive definite matrix. This algorithm is matrix-free and combines two existing diagonal estimators: based on the the Lanczos approach~\cite{chantas2008} and the Monte Carlo approach~\cite{bekas2007estimator}.

\subsection{Lanczos Diagonal Estimator} \label{ssec:lanczos}
Let $\ma \in \rnn$ be a symmetric positive definite (SPD) matrix. To estimate the diagonal of $\ma^{-1}$ we make use of a method introduced by \cite[Section IV]{chantas2008}, which uses the Lanczos iteration. We present this in a slightly different form using the preconditioned Lanczos iteration. Suppose we have a preconditioner  $\mg \in \rnn$ such that $\ma \approx \mg\mg^T$, $\mg$ is easy to invert, and the preconditioned matrix $\mg^{-1}\ma\mg^{-T}$ has a smaller condition number than $\ma$. Then after $k < n$ steps of the Lanczos iteration on the preconditioned matrix, assuming no breakdown, we have 
\[ \mg^{-1}\ma\mg^{-T}\mv_k = \mv_{k}\mt_k + \beta_{k+1}\vv_{k+1}\ve_k^T, \]
where $\mv_k\in \mathbb{R}^{n\times k}$ has orthonormal columns and $\mt_k \in \mathbb{R}^{k\times k}$ is a symmetric, tridiagonal matrix.  Since $\ma$ is SPD, $\mt_k$ is also SPD due to Cauchy interlacing theorem; therefore, we can compute the Cholesky factorization $\mt_k = \ml_k\ml_k^T$. Let us define $\mw_k = \mg^{-T}\mv_k\ml_k^{-T}$. After $n$ iterations of the Lanczos process, we have $\mv_n^T \mg^{-1}\ma\mg^{-T}\mv_n = \mt_n$, so that 
\[\ma^{-1} = \mg^{-T}\mv_n \mt_n^{-1} \mv_n^T\mg^{-1} = \mw_n\mw_n^T.\]
Suppose we define $\mb_k = \mw_k\mw_k^T$, then we can approximate $\ma^{-1} \approx \mb_k$. We can estimate the diagonals of $\ma^{-1}$ using the relation
\[\ve_i^T\ma^{-1} \ve_i \approx\ve_i^T\mb_k\ve_i =   \ve_i^T\mw_k\mw_k^T\ve_i = \sum_{j=1}^k w_{ij}^2. \]
Note that this is a monotonically increasing estimate for the $i-$th diagonal since $\mb_k = \mb_{k-1} + \vw_k\vw_k^T$. Note that we can easily adapt this algorithm to estimate the diagonals of $\mj\ma^{-1}\mj^T$ by instead working with $\mj\mb_k\mj^T$. In practice, the vectors which form the columns of $\mv_k$ lose orthogonality in finite precision, so we reorthogonalize them against previous vectors. A potential downside of this approach is that a large number of iterations $k$ may be required to estimate the diagonals accurately, even with the use of a preconditioner. Although this estimator does not require storing $\mw_k$, our approach in~\ref{ssec:new_approach} requires the storage of $\mw_k$ and this motivates the need to use a small $k$. 

\subsection{Monte Carlo Diagonal Estimator}\label{ssec:mc}
Another method to estimate $\diag(\ma^{-1})$ is using the matrix-free Monte Carlo diagonal estimator proposed in \cite{bekas2007estimator}. We describe this method to estimate the diagonals of $\ma^{-1}$, but this approach is applicable to any matrix $\mb$. In this approach, we draw a set of random vectors $\vz_\ell$ drawn from the Rademacher distribution, i.e., each vector has entries $\pm 1$ with equal probability. Then, we have the following approximation 
\begin{equation} \label{eq:MC_diagonal}
\vd^N_\text{MC} = \frac{1}{N}\left(\sum_{\ell=1}^N  \vy_\ell \odot \vz_\ell \right)  \qquad \vy_\ell= \ma^{-1}\vz_\ell, \quad 1 \leq \ell\leq N,
\end{equation}
where $\odot$ denotes elementwise multiplication. The vectors $\vy_\ell$  are computed in a matrix-free fashion with a CG solver applied to the system $\ma\vy_\ell = \vz_\ell$ for $1 \leq \ell \leq N$; solving these linear systems dominates the computational cost of the MC estimator. It can be shown that this estimate is unbiased in the sense that $\mathbb{E}[\vd^N_\text{MC}] = \diag(\ma^{-1})$, where $\mathbb{E}$ denotes the expectation. However, the convergence of this estimator is slow in the sense that a large number of vectors $N$ are required for an accurate estimate of the diagonals of $\ma^{-1}$. The authors in~\cite{bekas2007estimator} used a probing technique for estimating the diagonals of a matrix, that takes into account special properties of the matrix such as bandedness and decay of the magnitudes of the off-diagonal elements of the matrix. 
\subsection{LanczosMC Estimator}\label{ssec:new_approach}
Our approach for approximating $\diag(\ma^{-1})$ combines the ideas from Lanczos and Monte Carlo diagonal estimator. First, we approximate $\ma^{-1} \approx \mw_k\mw_k^T$ using the Lanczos method as described in~\ref{ssec:lanczos}. Then we decompose $\ma^{-1}$ as
\begin{equation}\label{eq:inv_formuala}
    \ma^{-1} = \mw_k\mw_k^T + \left(\ma^{-1} - \mw_k\mw_k^T\right).
\end{equation}
This equation implies that $\diag(\ma^{-1}) = \diag(\mw_k\mw_k^T) + \diag\left(\ma^{-1} - \mw_k\mw_k^T\right)$.
  Once the preconditioned Lanczos vectors are computed, the first summand can be calculated efficiently as it is a low-rank outer product. The second summand can be estimated via the MC method applied to $\ma^{-1} - \mw_k\mw_k^T$ rather than $\ma^{-1}$, as described in~\ref{ssec:mc}. The rationale behind this approach is that when $k$ is sufficiently large, we expect the approximation  $\ma^{-1} \approx \mw_k\mw_k^T$ to be reasonably accurate. However, numerical evidence suggests that the variance of the diagonal estimator applied to $\ma^{-1} - \mw_k\mw_k^T$  is much smaller than $\ma^{-1}$ resulting in more accurate estimators. Therefore, the use of Lanczos diagonal estimator can be interpreted as a variance reduction technique. This idea has been explored by~\cite{wu2016estimating,meyer2021hutch++} in the context of trace estimators. At this point, we lack a complete theoretical understanding of this but numerical experiments in~\ref{ssec:compex} provide convincing numerical evidence in support of our observation.  We summarize this new method for estimating $\diag(\ma^{-1})$ in Algorithm~\ref{alg:diag_new}.

  Related to the diagonal estimator, we can also estimate the trace of the matrix $\ma^{-1}$ as 
  \[ \text{trace}(\ma^{-1}) \approx \text{trace}(\mw_k^T\mw_k) + \frac{1}{N}\sum_{\ell=1}^N \vy_\ell^T\vz_\ell,\]
  where we have used the relation~\eqref{eq:inv_formuala} and the cyclic property of trace. 
\begin{algorithm}[!ht]                      
    \begin{algorithmic}[1]                    
        \REQUIRE  SPD matrix $\ma$, preconditioner $\mg$, integers $k$ denoting the number of iterations of the Lanczos approach and $N$ denoting the number of Monte Carlo samples.
        \RETURN An estimator $\vd^{k,N}_\text{LanczosMC}$ for $\diag(\ma^{-1})$.
            \STATE Run $k$ steps of the Lanczos approach on $\mg^{-1}\ma\mg^{-T}$ to compute $\mv_k$ and $\mt_k$. 
            \STATE Compute $\mw_k = \mg^{-T}\mb_k\ml_k^{-T}$ and 
                $$
                \vd^k_\text{Lanczos} = \diag(\mw_k\mw_k^T).
                $$
            \STATE Draw independent vectors $\{\vz_\ell\}_{\ell = 1}^N$ from the Radamacher distribution and compute $\vy_\ell = \left(\ma^{-1} - \mw_k\mw_k^T\right)\vz_\ell$ for $1 \leq \ell \leq N$. 
            \STATE Compute the Monte Carlo estimate
                $$
            \vd^N_\text{MC}= \frac{1}{N}\left(\sum_{\ell=1}^N  \vy_\ell \odot \vz_\ell \right) 
                $$
                where $\odot$ denotes elementwise multiplication.
            \STATE $\vd^{k,N}_\text{LanczosMC} = \vd^k_\text{Lanczos} + \vd^N_\text{MC}$
                \end{algorithmic}
    \caption{LanczosMC approach for estimating $\diag(\ma^{-1})$}
    \label{alg:diag_new}    
\end{algorithm}
    \paragraph{Computational Costs} To analyze the costs of Algorithm~\ref{alg:diag_new}, we make the following assumptions. We assume that the cost of applying the preconditioner $\mg^{-1}$ and $\mg^{-T}$ is the same and is denoted by $T_\text{prec}$. Furthermore, we denote the cost of applying $\ma$ and  $\ma^{-1}$ to a vector as $T_\text{apply}$  and $T_\text{solve}$ respectively. If an iterative method, such as CG, is used to apply $\ma^{-1}$, then the computational cost is $T_\text{solve} = N_\text{iter}(T_\text{apply} + 2T_\text{prec}) + \mathcal{O}(nN_\text{iter})$ flops; here $N_\text{iter}$ is the number of iterations taken by CG which may be different than $k$.    
    
    To calculate  $\vd^k_\text{Lanczos}$, we first perform $k$ steps of Lanczos on the preconditioned matrix $\mg^{-1}\ma\mg^{-T}$. The cost of running $k$ steps of this algorithm to generate $\mv_k$ and $\mt_k$ is $k(T_\text{apply}+ 2T_\text{prec}) + \mathcal{O}(nk)$ flops; the additional cost of computing $\mw_k$ and $\vd^k_\text{Lanczos}$ is $\mathcal{O}(nk^2)$ flops. Therefore, the total cost of computing $\vd^k_\text{Lanczos}$ is
    
    \[ T_\text{Lanczos} = k(T_\text{apply}+ 2T_\text{prec}) + \mathcal{O}(nk^2) \text{ flops}.\] 
    Each application of $\ma^{-1}-\mw_k\mw_k^T$ costs $T_\text{solve} + \mathcal{O}(nk^2)$, so the overall cost of computing $\vd^N_\text{MC}$ is 
    \[ T_\text{MC} = N(T_\text{solve} + \mathcal{O}(nk)) \text{ flops}.\]
    The total cost of the LanczosMC approach is 
    \[ \begin{aligned}T_\text{MC} =  & \>  T_\text{Lanczos} + T_\text{MC} \\
    =  &\> k(T_\text{apply}+ 2T_\text{prec}) + NT_\text{solve} + \mathcal{O}(nNk + nk^2) \text{ flops}. \end{aligned}\]
    In addition to the computational cost discussed here, the method requires storing $k$ vectors of length $n$. We briefly comment on the choice of parameters $k$ and $N$. The convergence of Monte Carlo estimators is slow, like $N^{-1/2}$, so we take $N$ to be large say $100-1000$. However, this computational cost can be trivially parallelized. The choice of $k$ depends on the storage budget and accuracy. In numerical experiments, we found that  $k\sim 100$ was a reasonable choice.

\subsection{Computational Examples}\label{ssec:compex}
To demonstrate the efficiency of \ref{ssec:new_approach} we estimate the diagonal of the inverse of 5 different test matrices. The test matrices were taken from \cite{suitesparse} and were chosen to represent a wide variety of sizes and condition numbers. Each matrix is sparse and SPD. We chose the matrices to be of relatively small size, so that we can study the accuracy of the estimators; however, our algorithms are applicable to larger problem sizes. In Table~\ref{tab:test_matrices} we list the name, size, and condition number for each test matrix.

\begin{figure}[H]
        \centering
         \includegraphics[scale=0.3]{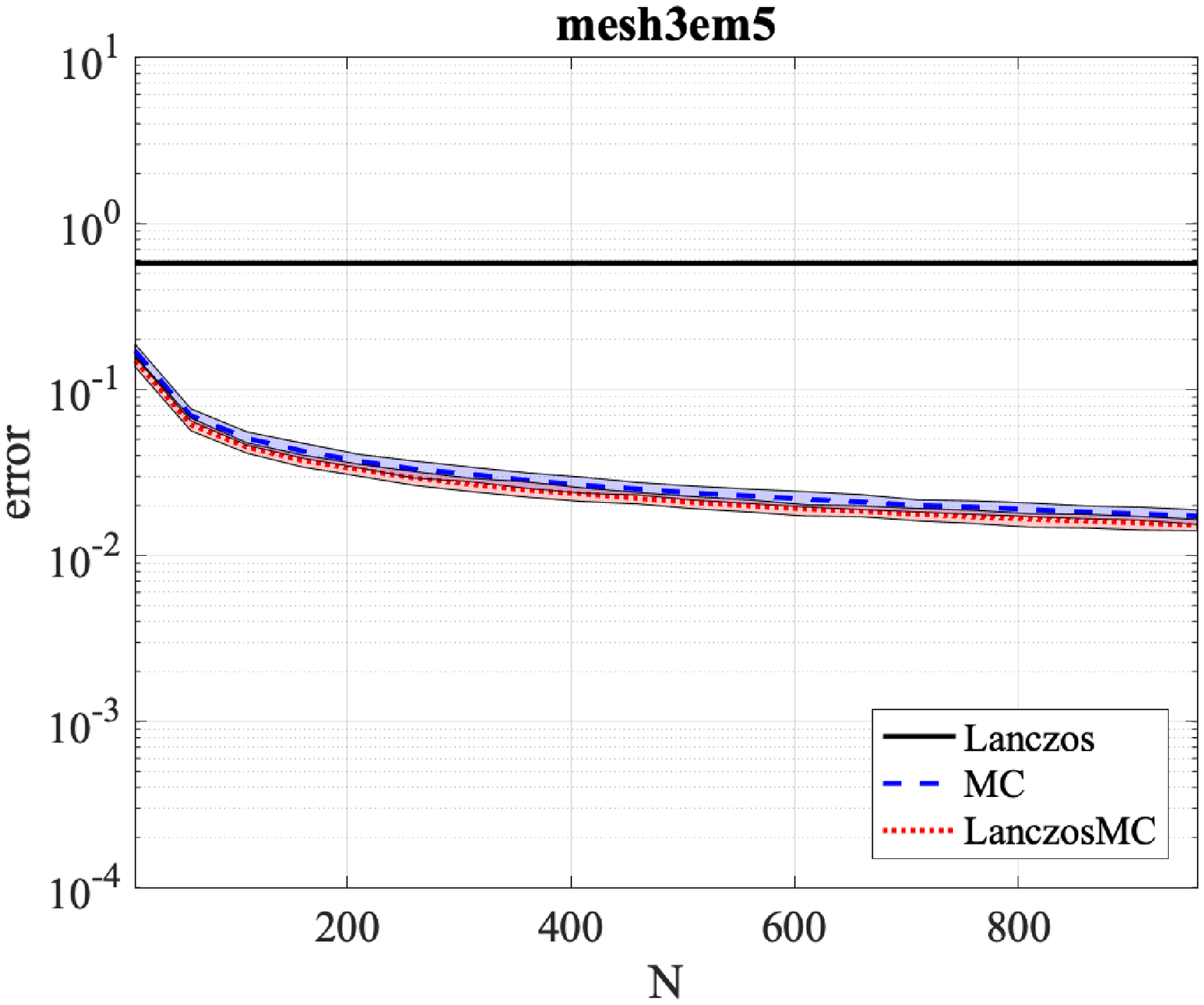} 
         \includegraphics[scale=0.3]{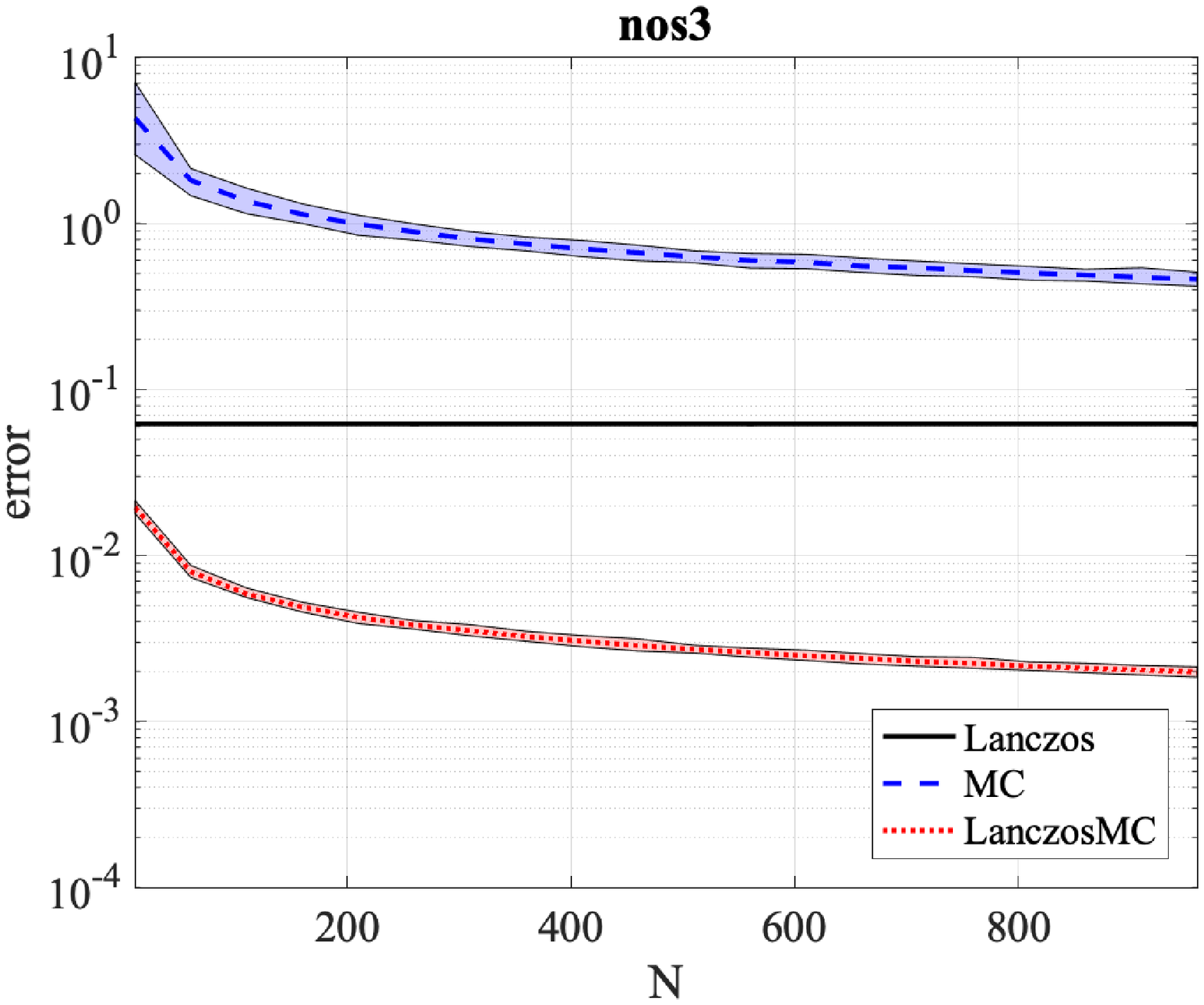} \\
        \includegraphics[scale=0.3]{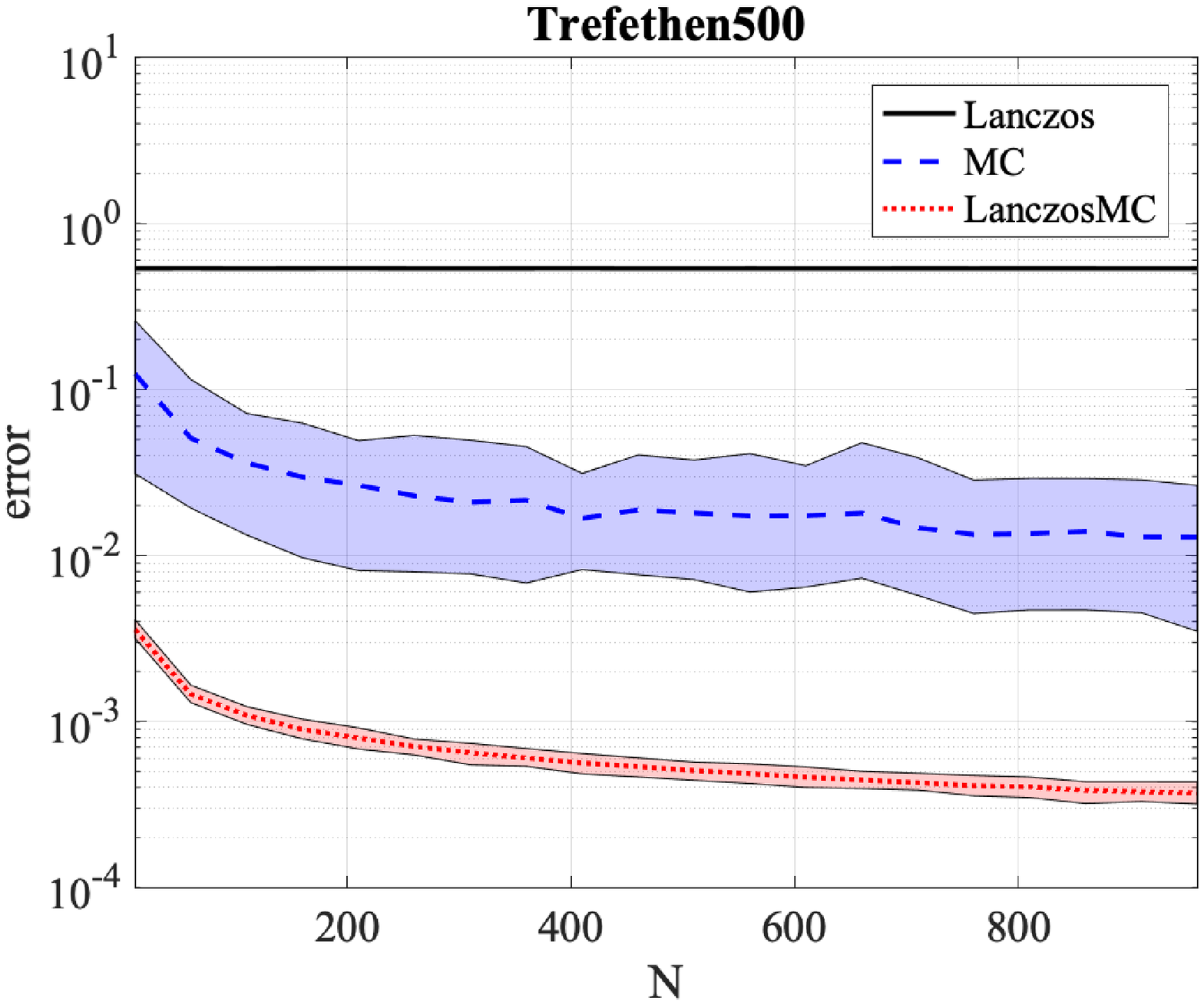} 
        \includegraphics[scale=0.3]{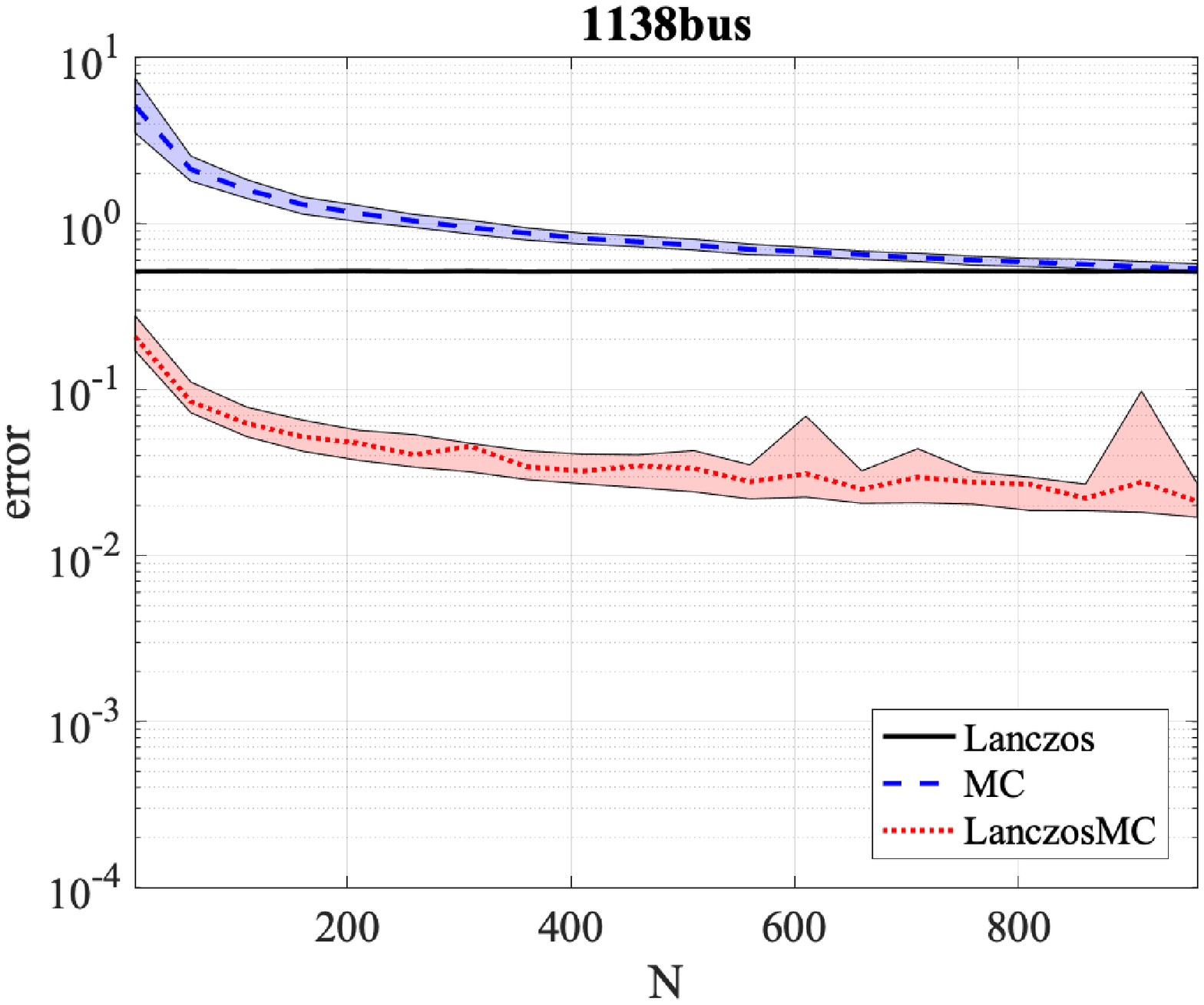} \\
        \includegraphics[scale=0.3]{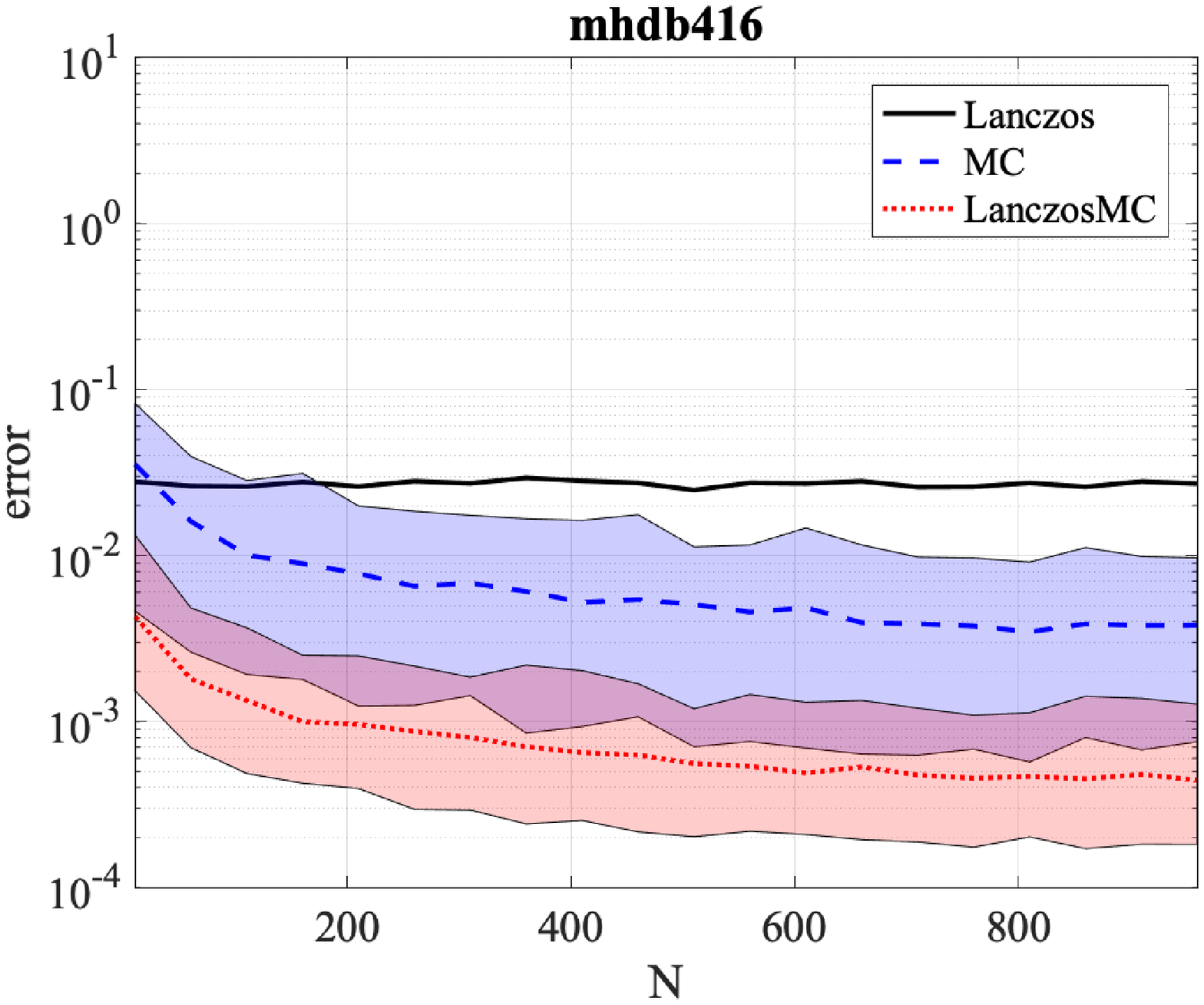}
        \caption{Relative error plots for the Monte Carlo, LanczosMC, and Lanczos estimators. The error statistics were calculated for 100 realizations for every fixed value of $N$.}
        \label{fig:diagestacc}
    \end{figure}

\paragraph{Experiment 1: Effect of varying $N$} We apply Algorithm~\ref{alg:diag_new} to the test matrices in Table~\ref{tab:test_matrices}. We chose $k=100$ and used the incomplete Cholesky preconditioner. We used CG to solve systems with $\ma$. We compute the relative error 
\[ \text{rel err} =  \frac{\norm{\vd_\text{C} - \vd_\text{True}}_2}{\norm{\vd_\text{True}}_2}
\]
where $\vd_\text{C}$ is the approximate diagonal and $\vd_\text{True}$ is the true diagonal.
The results are reported in Figure~\ref{fig:diagestacc}. The errors are averaged over $100$ different runs for a fixed sample size. In addition to plotting the average we also plot the $2.5$th and $97.5$th percentiles for the Monte Carlo and LanczosMC errors as well. From the plots it is clear that the LanczosMC estimator is accurate compared to MC estimator. It is also seen that the accuracy of the estimators improves as the condition number of the matrix increases. We see improvements in relative error by $1-2$ orders of magnitude (except the matrix mesh3em5, which is well-conditioned), suggesting that the Lanczos estimator is better for more ill-conditioned systems. %The slow convergence of the MC estimators suggests that 

    \begin{table}[H]
        \centering
        \begin{tabular}{c|c|c|c|c|c}
        %\hline
            \textbf{Matrix} & mesh3em5 & nos3 & {Trefethen500} & 1138bus & {mhdb416} \\
            \hline
            \textbf{Size} & $289 \times 289$ &   $960 \times 960$ & $500 \times 500$ &  $1138 \times 1138$ &  $416 \times 416$ \\
            \hline 
            \textbf{Condition Number} &   $5$ &   $4.63\times 10^{3}$ & $7.35 \times 10^{4}$  &  $1.23 \times 10^{7}$ & $5.052 \times 10^{9}$ \\
        \end{tabular}
        \caption{Summary of the test matrices used in the numerical experiments.}
        \label{tab:test_matrices}
    \end{table}
    
    In Figure \ref{fig:diagestacc} the shaded area between the $2.5$th and $97.5$th quantiles for the MC and Lanczos estimators corresponds to the 95\% confidence interval. For each test matrix we notice that the confidence interval for the LanczosMC estimator is smaller than the MC estimator for every test matrix. Since the confidence interval is a measure of variance in an estimator, this implies that the LanczosMC estimator has a smaller variance.
    
\begin{figure}[!ht]
        \centering
        \includegraphics[scale=0.35]{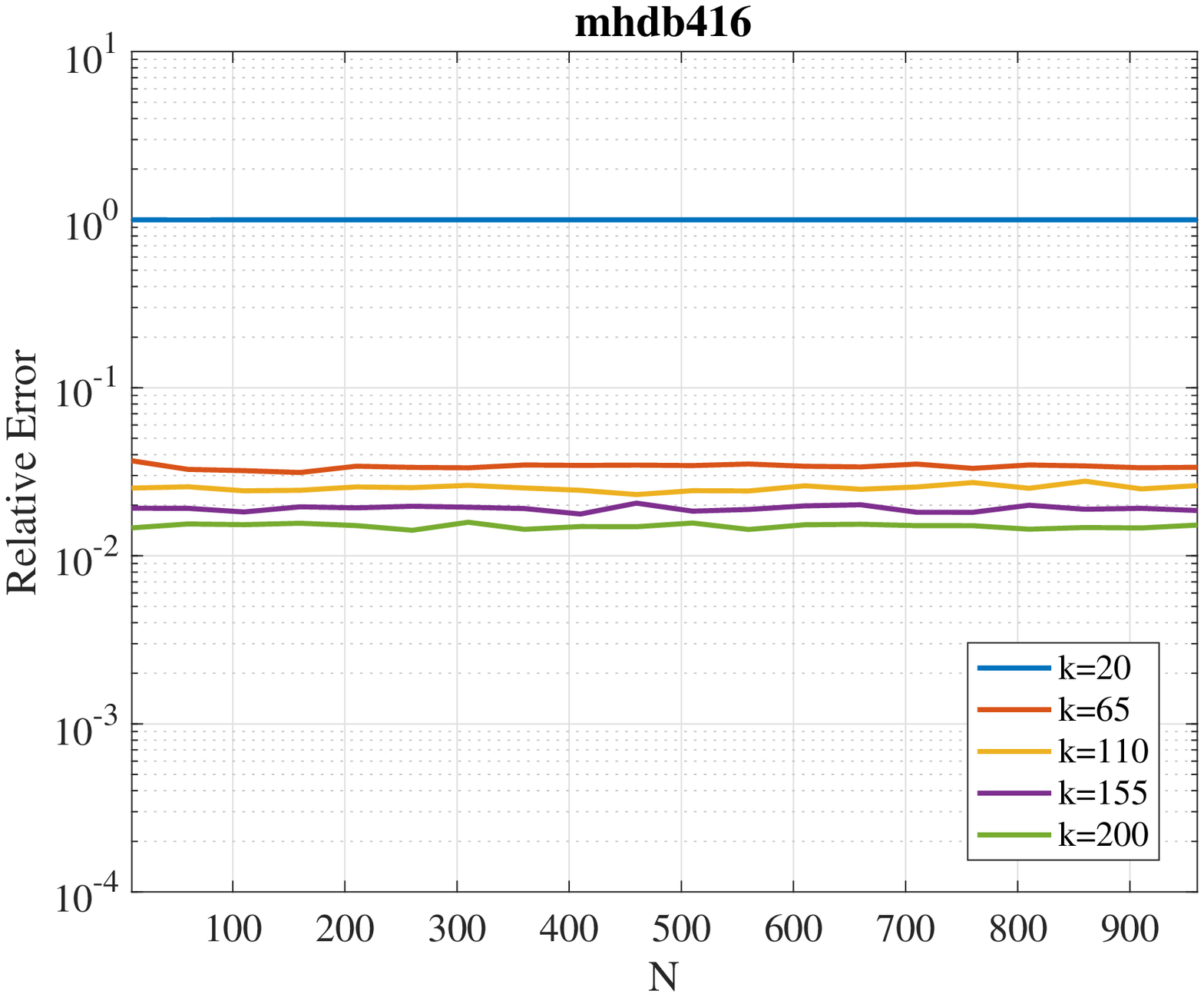}
        \includegraphics[scale=0.35]{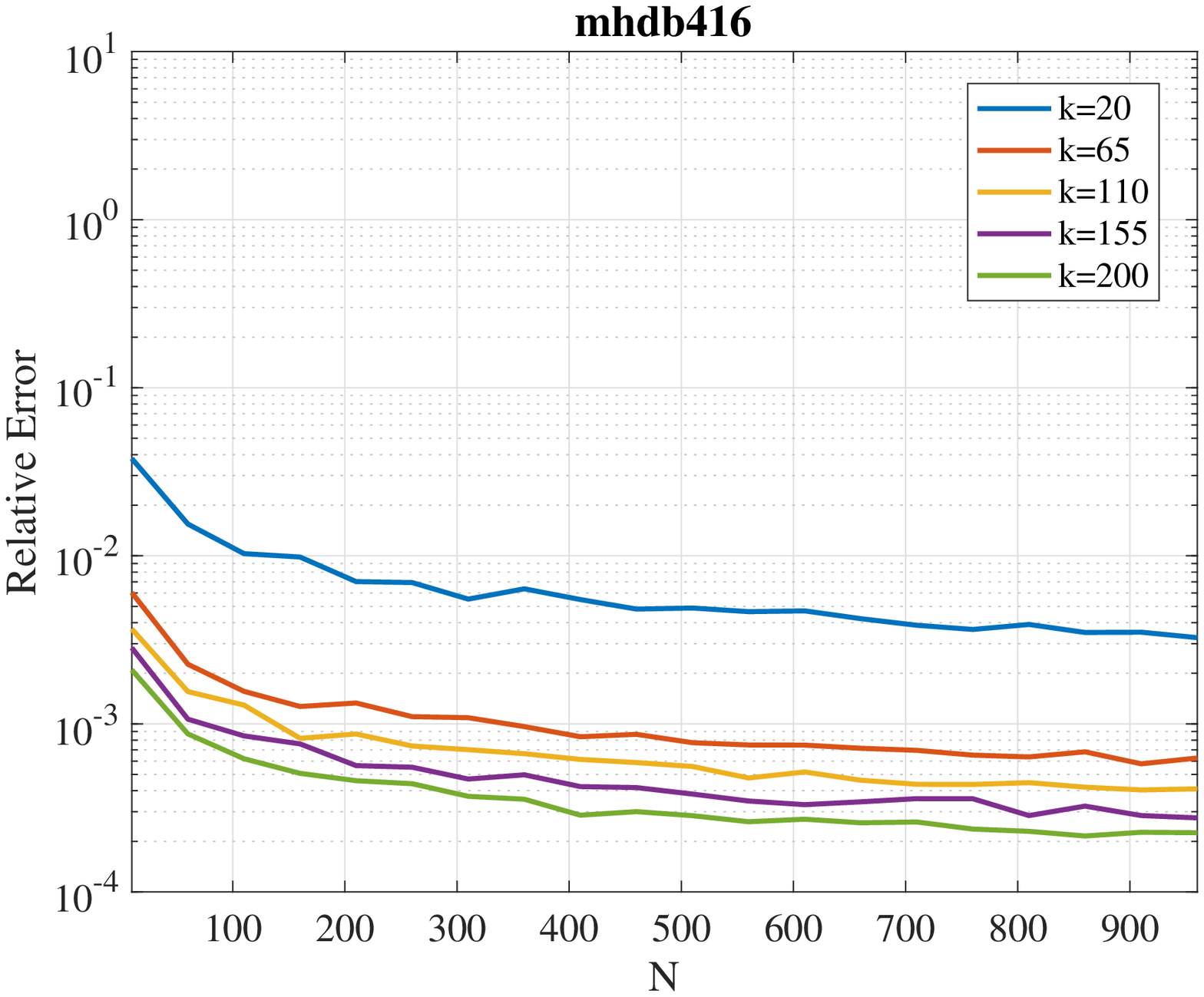}
        \caption{Relative error plots for different number of Lanczos iterations $k$: (left) mean of relative errors corresponding to Lanczos estimator, and  (right) Mean of relative errors corresponding to  LanczosMC estimator.}
        \label{fig:diagestak}
    \end{figure}

    \paragraph{Experiment 2: Effect of varying $k$} For this next experiment we analyze the impact of changing the number of Lanczos iterations $k$. We choose $\ma$ to be the mhdb416 matrix and plot the mean for the error while using the LanczosMC estimator. We also plot the mean of the error while using the Lanczos estimator. We choose the same preconditioner, range of Monte Carlo samples, and number of trials for each fixed number of Monte Carlo samples. The value of $k$ ranges from 20 to 200 and the results are reported in Figure \ref{fig:diagestak}.
    
\begin{comment}
As another mode of comparison we compute the error in approximating the diagonal using the Lanczos method. Similar to the previous experiment we fix $k$, plot the mean of the errors after 100 trials for a fixed, and then increase $N$. We do this for increasing values of $k$. The results are shown in the right panel of Figure \ref{fig:diagestak}. We see that as the number of Lanczos iterations increases the average relative error decreases. 
\end{comment}

\end{document}